\documentclass[final]{siamltex}
\usepackage{amssymb}
\usepackage{amsmath}
\usepackage{graphicx}
\usepackage{color}
\usepackage{graphicx, subfigure}
\usepackage[outer]{showlabels}
\newcommand{\rf}[1]{(\ref{#1})}

\newcommand{\tc}{\tilde{c}}
\newcommand{\tA}{\tilde{A}}

\newcommand{\tbt}{\tilde{b}^{T}}

\newcommand{\one}{\mathbf{1}}
\newcommand{\tb}{\tilde{b}}
\newcommand{\bt}{b^{T}}

\def\RR{\mathbb R}
\newcommand{\be}{\begin{equation}}
\newcommand{\ee}{\end{equation}}
\newcommand{\bea}{\begin{eqnarray}}
\newcommand{\eea}{\end{eqnarray}}
\newcommand{\bean}{\begin{eqnarray*}}
\newcommand{\eean}{\end{eqnarray*}}
\def\ba{\begin{array}{l}\displaystyle}
\def\ea{\end{array}}
\def\epsi{\varepsilon}

\definecolor{PBlue}{rgb}{0.0,0.0,0.69}
\definecolor{PGreen}{rgb}{0.0,0.69,0.0}
\definecolor{PMagenta}{rgb}{0.69,0.0,0.69}
\definecolor{PRed}{rgb}{0.8,0.0,0.0}
\newcommand{\acapo}{\\[-0.2cm]}
\title{Implicit-Explicit Runge-Kutta schemes for hyperbolic systems and kinetic equations in the diffusion limit}
\author{S. Boscarino\thanks{Mathematics and Computer Science Department, University of
Catania, Italy ({\tt boscarino@dmi.unict.it}).}\and L.
Pareschi\thanks{Mathematics Department, University of Ferrara,
Italy ({\tt lorenzo.pareschi@unife.it}).} \and
G.Russo\thanks{Mathematics and Computer Science Department,
University of Catania, Italy ({\tt russo@dmi.unict.it}).}}

\begin{document}
\maketitle
%\emph{ABSTRACT} \\
\begin{abstract}
We consider Implicit-Explicit (IMEX) Runge-Kutta (R-K) schemes for
hyperbolic systems with stiff relaxation in the so-called diffusion
limit. In such regime the system relaxes towards a convection-diffusion
equation. The first objective of the paper is to show that traditional partitioned IMEX R-K
schemes will relax to an explicit scheme for the limit equation with no need of modification
 of the original system. Of course the explicit scheme obtained in the limit suffers 
 from the classical parabolic stability restriction on the time step. 
 The main goal of the paper is to present an approach,
 based on IMEX R-K schemes, that in the diffusion limit relaxes to
 an IMEX R-K scheme for the convection-diffusion equation, in which 
 the diffusion is treated implicitly.
This is achieved by an original reformulation of the problem,
and subsequent application of IMEX R-K schemes to it. 
An analysis on such schemes to the reformulated problem shows that the
schemes reduce to IMEX R-K schemes for the limit equation, under the same
conditions derived for hyperbolic relaxation \cite{Bosc-10}.
Several numerical examples including neutron transport equations
confirm the theoretical analysis.
\end{abstract}

\noindent {\bf Key words.} IMEX Runge-Kutta methods, hyperbolic
conservation laws with sources, transport equations, diffusion
equations, stiff systems. \medskip

\noindent {\bf AMS subject classification.} 65C20, 65M06, 76D05,
82C40.

%\tableofcontents

\section{Introduction}
The development of numerical methods to solve hyperbolic systems
in diffusive regimes has been a very active area of research in
the last years (see for example \cite{JPT, Klar, LS, NP}).

Classical fields of applications involve diffusion in neutron
transport \cite{BC, DJ, GT, JPT2, LM}, drift-diffusion limit
in semiconductors \cite{JP, Klar3} and incompressible
Navier-Stokes limits in rarefied gas dynamic \cite{Klar2}. A
strictly related field of research concerns the construction of
schemes for the compressible Navier-Stokes limit (see \cite{BLM}
and the references therein). In such physical problems, the
scaling parameter (mean free path) may differ in several orders of
magnitude from the rarefied regimes to the diffusive regimes, and
it is desirable to develop a class of robust numerical schemes
that can work uniformly with respect to this parameter.

A prototype hyperbolic system of conservation laws with diffusive
relaxation that we will use to illustrate the subsequent theory is
the following, \cite{JPT, JP, NP} 
\begin{eqnarray} \label{I3} 
\begin{array}{l}
\displaystyle
u_t + v_x =0, \\[+.25cm]
\displaystyle v_t + \frac{1}{\epsi^2} p(u)_x = -\frac{1}{\epsi^2}
\left( v-q(u ) \right),
\\
\end{array}
\end{eqnarray}
where $p^\prime(u)>0$. System (\ref{I3}) is
hyperbolic with two distinct real characteristics speeds $\pm
\sqrt{p^\prime(u)}/\epsi$.

In the small relaxation limit, $\epsi\rightarrow 0$, the behavior
of the solution to (\ref{I3}) is, at least formally, governed by
the {\em convection-diffusion} equation
 \begin{eqnarray}  \label{I4}
\begin{array}{l}
\displaystyle u_t +
q(u)_x= p(u )_{xx},\\[+.25cm]
\displaystyle v= q(u ) -p(u)_x.
\\
\end{array}
\end{eqnarray} 
The so called subcharacteristic condition
\cite{CLL} for system (\ref{I3}) becomes
 \begin{eqnarray} | q^\prime (u ) |^2 < \frac{p^\prime(u)}{\epsi^2},
\label{I5} \end{eqnarray} 
and is naturally satisfied in the limit $\epsi\rightarrow 0$.

About the boundary conditions for system (\ref{I3}), we have to specify 
the domain in which we solve the problem. In a finite domain $x \in [a, b]$, 
one can use periodic boundary conditions, or can assign one independent 
condition at each boundary, since the two characteristic velocities have opposite
sign.

In practice, we shall assign two conditions at each boundary, one independent 
and the other compatible with the equations. Furthermore, we shall choose boundary conditions which
are compatible to the limit solution as $\varepsilon \to 0$. For example for system (\ref{I3}), if we 
set $v = q(u)$ at $x = a$ and $x=b$, compatibility with system (\ref{I3}) requires
\begin{eqnarray*}
\varepsilon^2 q'(u)v_x = p'(u)u_x, \ \ \ x = a,b.
\end{eqnarray*}
Such condition becomes $u_x = 0$ in the limit $\varepsilon \to 0$. In the sections on the numerical tests we shall specify the boundary conditions we use in each test.

In general, numerical approaches that work for hyperbolic system
with stiff relaxation terms do not apply directly in the diffusive
scaling since in these systems we have the presence of multiple
time-scales. 

In fact, together with the stiff relaxation term we
have a stiff convection term that contributes to the asymptotic
diffusive behavior. Then special care must be taken to ensure that
the schemes possess the correct zero-relaxation limit, in the
sense that the numerical scheme applied to system (\ref{I3})
should be a consistent and stable scheme for the limit system
(\ref{I4}) as the parameter $\varepsilon$ approaches zero independently of the discretization parameters. A
notion usually referred to as \emph{asymptotic preservation}.
For a nice survey on asymptotic preserving scheme for various kinds of
systems see, for example, the review paper by Shi Jin (\cite{SJAP}).
Furthermore, a different approach to the derivation of asymptotic preserving schemes is described in the review by Pierre Degond \cite{DeG}.
In the case of Boltzmann kinetic equations we also refer to the recent review by two of the authors \cite{PR_VKI}.

IMEX Runge-Kutta (R-K) schemes \cite{ARS, Bosc-07, Bosc-09, CK, A}
represent a powerful tool for the time discretization of such
stiff systems. Unfortunately, since the characteristic speed of
the hyperbolic part is of order $1/\varepsilon$, standard IMEX R-K
schemes developed for hyperbolic systems with stiff relaxation
\cite{A, Bosc-10} become useless in such parabolic scaling,
because the CFL condition would require $\Delta
t=\mathcal{O}(\varepsilon \Delta x)$. Of course, in the diffusive
regime where $\varepsilon < \Delta x$, this is very restrictive
since for an explicit method a parabolic condition $\Delta t
=\mathcal{O}(\Delta x^2)$ should suffice.

Most previous works on asymptotic preserving schemes for
hyperbolic systems and kinetic equations with diffusive relaxation
focus on schemes which in the limit of infinite stiffness become
consistent explicit schemes for the diffusive limit equation
\cite{JPT, Klar, LS, LM, NP}. Such schemes have been derived by splitting
the stiff hyperbolic part into an explicit (non-stiff) term, and an implicit (stiff)
term. Here we show that by applying partitioned IMEX R-K schemes, in which 
 the stiffness is associated with the variable and not with the operator,
 one obtains IMEX R-K schemes that naturally relax to the explicit
 scheme applied to the limit convection-diffusion equation. All these explicit schemes clearly suffer from the usual stability restriction $\Delta t =\mathcal{O}(\Delta
x^2)$.

In this paper we present a general methodology to overcome such stability restriction 
which applies to a broad class of problems. The idea
is to reformulate problem \rf{I3} by properly combining the
limiting diffusion flux with the convective flux. This allows to
construct a class of IMEX R-K schemes that work with high order
accuracy in time and that, in the diffusion limit (i.e.\ when
$\varepsilon\to 0$), originate an IMEX method for the limiting
convection-diffusion equation \rf{I4}. Other reformulations whose
 goal is to obtain asymptotic-preserving methods have been
proposed in \cite{JPT, JF}. Schemes that avoid such time step restriction and provide
fully implicit solvers in the case of transport equations have
been analyzed in \cite{BC}. 

Our new approach allows a hyperbolic CFL condition
$\Delta t = \mathcal{O}(\Delta x)$ independent of $\varepsilon$
when applied to \rf{I3} in all regimes. The aim of this paper is
to derive and analyze different types of IMEX R-K schemes when
applied to the reformulated problem in the stiff regime
($\varepsilon \rightarrow 0$).

The rest of the paper is organized as follows. The next section is devoted to partitioned  IMEX R-K schemes. It is shown 
that they relax to the explicit scheme applied to the convection diffusion limit. 
In Section \ref{Sec3} the new approach is introduced and analyzed.
In particular, following \cite{Bosc-10}, we prove that under suitable assumptions the
IMEX R-K schemes are consistent with the diffusion limit. The analysis is based
on a power series expansion in $\varepsilon$ of the exact and numerical solution.
It is shown that, at lowest order in $\varepsilon$, the model system is a set of differential-algebraic
equations of index 1, i.e. it can be transformed into a set of ordinary differential equations
by one time differentiation. Compatibility between exact and numerical solution at different orders
in $\varepsilon$ introduces additional order conditions on the coefficients of the IMEX schemes.
To the lowest order such conditions are referred as \emph{index 1 order conditions}.

After a short section on space discretization obtained by conservative
 finite difference schemes, in Section \ref{Sec5} we report several numerical examples and tests. In Section
  \ref{Sec6} we consider one-dimensional neutron transport equation and present several
numerical results and comparison with schemes available in the literature. Additional
technical material is given in the separate Appendices.
%%%%%%%%%%%%%

\section{Partitioned IMEX R-K schemes} \label{Sec2}
The first observation is that in system (\ref{I3}) the stiffness is naturally associated to the variable $v$ rather then to some operator.
The system has the structure of a singular perturbation problem \cite{HW}, and it can be treated by a partitioned R-K scheme
 in which the first equation is treated explicitly and the second implicitly
\begin{eqnarray} \label{I33}
\begin{array}{lr}
\displaystyle
u_t  = - v_x, &\textrm{(Explicit)}\\[+.25cm]
\epsi^2\displaystyle v_t =-( p(u)_x + v-q(u )). &\textrm{(Implicit)}
\\
\end{array}
\end{eqnarray}
This approach has been used, for example, in \cite{Klar, NP2}.

By using a method of lines approach (MOL), we discretize system (\ref{I33}) in space by a uniform mesh $\left\{x_i\right\}_{i  =1}^N$  and  $U_i(t) \approx u(x_i, t)$, $V_i(t) \approx v(x_i, t)$. We obtain a large system of ODE's 
\begin{eqnarray}\label{MLa}
\begin{array}{l}
\displaystyle
U_t  = \mathcal{F}(V), \\[+.25cm]
\epsi^2\displaystyle V_t = G(U) - V,
\\
\end{array}
\end{eqnarray}
with  $U(t) = (U_1(t), U_2(t), ..., U_N(t))^T \in \mathbb{R}^N$ and $V(t) = (V_1(t), V_2(t), ..., V_N(t))^T \in \mathbb{R}^N$, 
where $\mathcal{F}(V) = -D V$ and $G(U)= Q(U)-D p(U) $. Here $D V$ and $D p(U)$ (with a slight abuse of notation) denote  
the discretization of the convective terms $v_x$, $p(u)_x$, while $Q(U)$ represents the discretization of the term $q(u)$.

As we shall see, if the implicit scheme is $L$-stable \cite{HW}, in the 
limit $\varepsilon \rightarrow 0$ the IMEX R-K scheme will relaxe to the explicit scheme applied to the limit equation
\begin{eqnarray}\label{MLalimit}
U_t =: \hat{\mathcal{F}}(U),%G(U)).
\end{eqnarray}
where $\hat{\mathcal{F}}(U) = \mathcal{F}(G(U))$.

For example, implicit-explicit Euler scheme applied to system (\ref{MLa}) gives
\begin{eqnarray*}
U_{n+1} &=& U_{n} + \Delta t \mathcal{F}(V_{n})\\
V_{n+1} &=& \frac{\varepsilon^2 V_n + \Delta t G(U_{n+1})}{\varepsilon^2+ \Delta t},
\end{eqnarray*}
where we have discretized the interval of integration by a time mesh $\left\{t_n\right\}_{n  =1}^\mathcal{N}$ and $U_n \approx  U(t_n)$.
As $\varepsilon \rightarrow 0$, $V_{n+1} = G(U_{n+1})$ and therefore $U_{n+1} = U_{n}+\Delta t \hat{\mathcal{F}}(U_{n})$.

In the case $p(u) = u$ and $q(u)=0$, the method would relax to explicit Euler scheme applied to the diffusion equation, thus suffering the usual parabolic 
CFL stability condition $\Delta t \leqslant \Delta x^2/2$.
This approach will be denoted as \emph{partitioned approach}. 

\subsection{Classification of IMEX R-K schemes}
IMEX R-K schemes have been widely used in the literature
to treat problems that contain both stiff and non stiff terms
\cite{ARS,CK,A}. The stiff terms are treated implicitly, while the
non stiff terms are treated explicitly, thus lowering the
computational complexity of the scheme.

Usually such a scheme is characterized by the $s \times s$ matrices $\tilde{A} = (\tilde{a}_{ij})$, $A
= ({a}_{ij})$ and the vectors $\tilde{b}$, $b$ $\in \RR^s$, and can be represented by a double
$tableau$ in the usual Butcher notation
\begin{displaymath}
\begin{array}{c|c}
\tilde{c} & \tilde{A}\\
\hline
\\[-.3cm]
 & \ \ \tilde{b^T} \end{array} ,\quad
\begin{array}{c|c}
{c} & {A}\\
\hline
\\[-.3cm]
 & \ \ {b^T} \end{array}.
\end{displaymath}
The coefficients $\tilde{c}$ and $c$ are used if the right hand side depends explicitly on time.
We assume that they satisfy the usual relation
\begin{eqnarray}\label{eq:candc}
\tilde{c}_i = \sum_{j=1}^{i-1} \tilde a_{ij}, \quad c_i = \sum_{j=1}^{i} a_{ij}.
\end{eqnarray}
Matrix $\tilde{A}$ is lower triangular with zero diagonal, while
matrix $A$ is lower triangular, i.e. the implicit scheme is a
diagonally implicit Runge-Kutta (DIRK). This choice guarantees
that the term $\mathcal{F}(V)$ in (\ref{MLa}) is always explicitly evaluated.

IMEX R-K schemes presented in the literature can be classified in
two main different types characterized by the structure of the
matrix $A = (a_{ij})_{i,j=1}^s$ of the implicit scheme.
\begin{definition}
\label{Def:A}
We call an IMEX R-K method of type A (see \cite{A}) if the matrix $A \in \RR^{s \times s}$ is invertible.\\
\end{definition}
\begin{definition}
\label{Def:CK} We call an IMEX R-K method of type CK (see
\cite{CK}) if the matrix $A \in R^{s \ \times \ s}$ can be written
as
\begin{eqnarray*}
A = \left(\begin{array}{ll} 0 & 0\\
                        a  &  \hat{A}\end{array}\right)
               \end{eqnarray*}
with $a\in \RR^{(s-1)}$ and the submatrix $\hat{A} \in \RR^{(s-1)
\ \times \ (s-1)}$ invertible. In the special case $a=0$ the
scheme is said to be of type ARS (see \cite{ARS}).
\end{definition}

We note that schemes CK are very attractive because they allow some simplifying assumptions, 
that make order conditions easier to treat, therefore permitting the construction of higher order
IMEX R-K schemes. On the other hand, schemes of type A are more amenable to a theoretical 
analysis, since the matrix $A$ of the implicit scheme is invertible. This is why we start our analysis
with the latter schemes.  
\subsection{Analysis of IMEX schemes for the partitioned approach}
Now as an example we perform the analysis of type A scheme when applied to system (\ref{I33}). A 
similar analysis is possible also for CK schemes. We will restrict our analysis to the limit case $\varepsilon \to 0$. 

Applying an IMEX R-K scheme to system (\ref{MLa}) we obtain
\begin{eqnarray}\label{typeAPart}
  U_{n+1} & = & U_n + \Delta t \sum_{k=1}^{s} \tilde{b}_{k} \mathcal{F}(\mathcal{V}_k),
   \nonumber         \acapo
            &   &                     \acapo %\label{eq:IMEXF4}
  \varepsilon^2   V_{n+1}  & = &    \varepsilon^2 V_n +  \Delta t \sum_{k=1}^{s} b_k  (G(\mathcal{U}_k)-\mathcal{V}_k),      \nonumber
\end{eqnarray}
for the numerical solution and
\begin{eqnarray}\label{stypeAPart}
    \mathcal{U}_k & = &   U_n + \Delta t\sum_{j=1}^{k-1} \tilde{a}_{kj} \mathcal{F}(\mathcal{V}_j)   \nonumber         \acapo
            &   &                     \acapo %\label{eq:IMEXF4bis}
  \varepsilon^2   \mathcal{V}_k  & = & \varepsilon^2  V_n + \Delta t \sum_{j=1}^{k} a_{kj} (g(\mathcal{U}_j)-\mathcal{V}_j),      \nonumber
\end{eqnarray}
for the internal stages.

By Definition \ref{Def:A} and $A$ invertible we obtain from the second equation in (\ref{typeAPart})
\begin{eqnarray}\label{gPar}
\Delta t (G(\mathcal{U}_k) - \mathcal{V}_k) = \varepsilon^2\sum_{j=1}^{k}\omega_{kj}(\mathcal{V}_j-V_n).
\end{eqnarray}
From now on, $\omega_{kj}$ are
the elements of the inverse matrix $A^{-1}$. Now inserting (\ref{gPar}) into the numerical solution $V_{n+1}$ and setting $\varepsilon = 0$, we get
\begin{eqnarray}%\label{index1APar}
     V_{n+1}  & = & \mathcal{R}(\infty)  V_n +  \Delta t \sum_{k=1}^{s} b_k \omega_{kj}  \mathcal{V}_j.
\end{eqnarray}
Here we denoted by
\[\mathcal
{R}(\infty) = 1-b^TA^{-1}\one = \lim_{z\rightarrow \infty} \mathcal{R}(z),\]
where $\mathcal{R}(z)$ is \emph{the stability function} of the implicit
scheme defined by (see \cite{HW}, Sect. IV.3)
\begin{eqnarray}\label{Sfunc}
\mathcal{R}(z) = 1 +
zb^{T}(I-zA)^{-1}{\bf 1},
\end{eqnarray}
with $b^{T} = (b_1,\ldots,b_s)$ 
and $\one = (1,\ldots,1)^{T}$.\\

By (\ref{gPar}) we get $\mathcal{V}_k = G(\mathcal{U}_k)$ when $\varepsilon = 0$. This yields that $\hat{\mathcal{F}}(\mathcal{U}_k) = \mathcal{F}( G(\mathcal{U}_k))$, and we obtain
\begin{eqnarray*}
U_{n+1} & = &  U_n  +  \Delta t \sum_{k = 1}^{s} \tilde{b}_k  \hat{\mathcal{F}}(\mathcal{U}_k) \nonumber         \acapo
\end{eqnarray*}
with
\begin{eqnarray*}          
              \mathcal{U}_k & = &  U_n + \Delta t \sum_{j = 1}^{k-1} \tilde{a}_{kj} \hat{\mathcal{F}}(\mathcal{U}_j)    \nonumber         \acapo
\end{eqnarray*}
internal stages and $k = 1,\ldots,s$. 
This represents the explicit scheme of the starting IMEX R-K one of type A applied to the limit equation (\ref{MLalimit}) obtained by (\ref{MLa}) when $\varepsilon \to 0$. As particular case, if $p(u) = u$ and $q(u)=0$, this is the explicit scheme applied to the limit diffusion equation under the usual parabolic stability restriction ($\Delta t \leqslant \Delta x^2/2$).

\section{Overcoming parabolic stiffness}\label{Sec3}
In order to overcome such stability restriction, we reformulate system (\ref{I3}) as the
equivalent system \begin{eqnarray}
\begin{array}{l}
\displaystyle
u_t  =- (v+ \mu p(u)_x)_x + \mu p(u)_{xx}, \\[+.25cm]
\epsi^2\displaystyle v_t =- p(u)_x -v+q(u ),
\end{array}
\label{I3b} \end{eqnarray} 
where the term $\mu p(u)_{xx}$ has been added and
subtracted to the first equation in (\ref{I3}). Here
$\mu=\mu(\varepsilon)\in [0,1]$ is a free parameter such that
$\mu(0)=1$. The idea is that, since the quantity $v+p(u)_x$ is
close to $q(u)$ as $\varepsilon\to 0$, the first term on the right
hand side can be treated explicitly in the first equation, while
the term $p(u)_{xx}$ will be treated implicitly. This can be done
naturally by using an Implicit-Explicit approach, as we will
explain later. Let us point out that the choice $\mu\equiv 1$, as
shown in Appendix \ref{A1} for a first order implicit-explicit scheme, guarantees the largest stability region of
the method.

Next we will study the behavior of the different IMEX R-K schemes
when applied to system (\ref{I3b}) in the diffusion
limit. In particular we will show that such schemes reduce to the same IMEX R-K schemes for the limit equation and
no parabolic stability
restriction on the time step appears in the diffusive limit.

%%%%%%%%%%%%%%%%%%%%%%%%%%%%%%%%%%%%%
\subsection{The new approach}
System (\ref{I3b}) can be written in the form
\begin{eqnarray}\label{SPP}
\begin{array}{l}
\displaystyle u' = f_1(u,v) + f_2(u), \\ [+.25cm]
\displaystyle {\varepsilon^2}v' =  g(u,v)\\
\end{array}
\end{eqnarray}
where the primes denote the time derivatives and  \[f_1(u,v) =
-( v + \mu p(u)_x)_x,\quad f_2(u) =
\mu p(u)_{xx},\]
\[g(u,v) =- p(u)_x - v+q(u).\] 
Notice that, throughout this paper, $g(u,v)$ (and therefore $g_v(u,v)$), depends only
algebraically on $v$, while it may contain differential operators acting on $u$.

Now we apply an IMEX-RK scheme to
system (\ref{SPP})
 where $(f_1(u,v),0)^T$ is evaluated explicitly and $(f_2(u),g(u,v))^T$ implicitly.
Note that if $f_2(u)$ is evaluated explicitly then by cancelation
the IMEX-RK scheme will reduce to the typology of asymptotic
preserving methods studied in \cite{BPR, NP2}.

In the limit $\varepsilon \rightarrow 0$ from (\ref{SPP}) we obtain a differential
algebraic system (DAE)
\begin{eqnarray}\label{DA1}
\begin{array}{l}
\displaystyle u' = f_1(u,v) + f_2(u), \\ [+.25cm]
\displaystyle 0 = g(u,v). \\
\end{array}
\end{eqnarray}
In order to guarantee the solvability of system (\ref{DA1}) we
assume that the Jacobian matrix $g_v(u,v)$ is invertible, and then
the DAE is said to be of index one \footnote{The index of a DAE is
the number of times one has to differentiate the function $g$ to
obtain a system of ODE's. For example, differentiating the
function $g$, one obtains $g_u(u,v)u' + g_v(u,v)v' = 0$. If $g_v$
is invertible, system (\ref{DA1}) can be written as $u' = f(u,v)$,
$v' = -g_v'g_u f$.}. Note that if $g_v$ has a bounded inverse in a
neighborhood of the exact solution, we can use the inverse
function theorem in order to write \[v(t) = G(u(t))\] for some
$G(u)$ which inserted into $u' = f_1(u,v) + f_2(u)$ gives $u' = f_1(u,G(u)) + f_2(u)$.
From now on we always assume that this is the case.
Then, as $\varepsilon\to 0$ system (\ref{SPP}) reduces to
\begin{eqnarray}\label{SPP3}
u' &=& \hat{f}_1(u) + f_2(u),
\end{eqnarray}
where $\hat{f}_1(u) = f_1(u,G(u))$ and $v = G(u)$.
This approach will be denoted \emph{BPR approach}.

First order implicit-explicit Euler scheme that uses this approach
is reported in Appendix \ref{A1}, where a stability analysis is performed.
In particular it is shown that as $\varepsilon \to 0$, the parabolic restriction on time 
step disappears. 

In the sequel we restrict our analysis to the limit case
$\varepsilon \rightarrow 0$ where the main goal is to capture the
diffusive limit.

\subsection{Analysis of TYPE A IMEX schemes}
Applying an IMEX R-K scheme of type A to system (\ref{SPP}) we obtain
\begin{eqnarray}\label{typeA}
  u_{n+1} & = & u_n + \Delta t \sum_{k=1}^{s} \tilde{b}_{k} f_1(U_k,V_k) + \Delta t \sum_{k = 1}^{s} b_k f_2(U_k)
   \nonumber         \acapo
            &   &                     \acapo %\label{eq:IMEXF4}
  \varepsilon^2   v_{n+1}  & = &    \varepsilon^2 v_n +  \Delta t \sum_{k=1}^{s} b_k  g(U_k,V_k),      \nonumber
\end{eqnarray}
for the numerical solution and
\begin{eqnarray}\label{stypeA}
    U_k & = &   u_n + \Delta t\sum_{j=1}^{k-1} \tilde{a}_{kj} f_1(U_j,V_j) +
  \Delta t \sum_{j = 1}^{k} a_{kj}  f_2(U_j)   \nonumber         \acapo
            &   &                     \acapo %\label{eq:IMEXF4bis}
  \varepsilon^2   V_k  & = & \varepsilon^2  v_n + \Delta t \sum_{j=1}^{k} a_{kj} g(U_j,V_j),      \nonumber
\end{eqnarray}
for the internal stages (notice a slight changes of notation with respect to Section \ref{Sec2}).

Starting from (\ref{typeA}) and (\ref{stypeA}), by Definition (\ref{Def:A}) and $A$ invertible, we obtain from the second equation in (\ref{stypeA})
\begin{eqnarray*}
\Delta t g(U_k, V_k) = \varepsilon^2\sum_{j=1}^{k}\omega_{kj}(V_j-v_n),
\end{eqnarray*}
Inserting this into the numerical solution $v_{n+1}$ we make $v_{n+1}$ independent of $\varepsilon^2$ and setting $\varepsilon = 0$, we get
\begin{eqnarray}\label{index1A}
  u_{n+1} & = &  u_n  +  \Delta t \sum_{k = 1}^{s} \tilde{b}_k \hat f_1(U_k) +  \Delta t \sum_{k = 1}^{s} b_k f_2(U_k) \nonumber         \acapo
            &   &  \acapo
            &   &                    \nonumber \acapo
   v_{n+1}  & = & \mathcal{R}(\infty)  v_n +  \Delta t \sum_{k=1}^{s} b_k \omega_{kj}  V_j,      \nonumber
\end{eqnarray}
with $\hat f_1(U_k) = f_1(U_k, G(U_k))$, and stage values 
\begin{eqnarray}\label{index1bisA}
   U_k & = &  u_n + \Delta t \sum_{j = 1}^{k-1} \tilde{a}_{kj} \hat f_1(U_j) + \Delta t \sum_{j = 1}^{k} a_{kj} f_2(U_j)   \nonumber         \acapo
            &   &                    \acapo
  0 & = & g(U_k, V_k).     \nonumber
\end{eqnarray}
for $k = 1,\ldots,s$. The latter equality implies $V_k=G(U_k)$, $k=1,\ldots,s$.

Note that if we require that the implicit part of the scheme is
{\em stiffly accurate}, i.e. if \[ b^TA^{-1} = e_s^T,\] where $e_s =
(0,\ldots,0,1)^T$, then by (\ref{Sfunc}) \[\mathcal{R}(\infty) =  1 -\bt A^{-1}{\bf 1} =
1-e_s^T{\bf 1} = 1-1 = 0. \] This implies that if the implicit
scheme is $A$-stable and stiffly accurate it is also $L$-stable and
$ v_{n+1}  =    V_s = G(U_s)$.

It is interesting to note that, if we consider system (\ref{I3b})
with $q(u) = 0$, when $\varepsilon = 0$ we get
a purely diffusive system which means that the term $f_1(u,v)$ in
(\ref{SPP}) disappears. Therefore, by BPR approach, the IMEX R-K scheme of type A in the limit $\varepsilon \rightarrow 0$ 
becomes a stiffly accurate DIRK scheme and hence no stability restriction on the time step is required in the diffusive
limit, i.e. we got an unconditionally stable method. Another advantage of this new approach is the following. Usually 
the numerical solution $(u_{n+1},v_{n+1})$ in (\ref{index1A}) in the case $\varepsilon = 0$ will not lie on the
manifold $g(u,v) = 0$ since $g(u_{n+1},v_{n+1})$ is not necessarily zero. But this approach 
guarantees that in the limit $\varepsilon \to 0$ we obtain a stiffly accurate implicit scheme  
and hence $u_{n+1} = U_s$, implying $g(u_{n+1},v_{n+1}) = 0$.

In the general case of systems for which $q(u) \neq 0$, it is $f_1(u,v) \neq
0$ and, by using the BPR approach, in the limit case $\varepsilon \rightarrow 0$ we obtain an IMEX R-K scheme 
with a non vanishing explicit term in which the diffusion term $f_2(u)$ is treated implicitly and a classical CFL hyperbolic 
condition for the time step is required. In general $g(u_{n+1}, v_{n+1}) \neq 0$ even if all stage values lie on the
manifold, (see the second equation in (\ref{index1bisA})).
However, if the explicit scheme has the property that
$u_{n+1} = U_s$, and the implicit scheme is stiffly accurate,
then, in the limit  as $\varepsilon \rightarrow 0$ the numerical
solutions are projected on the manifold $g(u_{n+1},v_{n+1}) = 0$, because
$g(u_{n+1},v_{n+1}) = g(U_s,V_s) = 0$.

From the above discussion it is clear that the property $u_{n+1} =
U_s$ is crucial if we want that the numerical solution is
projected to the limit manifold as $\varepsilon \rightarrow 0$.
We emphasize that  there is a class of
$s$-stage explicit R-K methods for which $u_{n+1} =
U_s$; such methods are called {\em First Same As Last}
(FSAL), and their coefficients satisfy $a_{s,i} = b_i$ for $i =
1,\ldots,s-1$ and $b_s = 0$.  They have the advantage of requiring $s-1$
function evaluations for each step (see  \cite{HNW} for details). 
FSAL methods are often used in the contest of embedded methods, such as  
the popular Dormond-Prince method (DOPRI) \cite{DP},
on which MATLAB routine ode45 is based on.

From the arguments above, in order to capture the limit as $\varepsilon \rightarrow 0$, 
it is important that the implicit part on an IMEX R-K is stiffly accurate  and the explicit part is FSAL.
This motivates the following 

\begin{definition}\label{1def} We say that a IMEX R-K scheme
is \emph{globally stiffly accurate} if $\bt = e_s^T A$ and $\tbt =
e_s^T\tA$, with $e_s = (0,\ldots,0,1)^T$, and $c_s = \tilde{c}_s =
1$, i.e. the numerical solution is identical to the last internal
stage value of the scheme.
\end{definition}

From (\ref{typeA}) and (\ref{stypeA}) we observe that if an IMEX
R-K is globally stiffly accurate, then
$u_{n+1} = U_s$, $v_{n+1} = V_s$, and therefore $\lim_{\varepsilon
\rightarrow 0} g(u_{n+1}, v_{n+1}) = 0$. 

\paragraph{General remarks for type A}

\begin{itemize}
\item It is worth mentioning some important aspects
about type A schemes. First of all, in \cite{Bosc-07} Boscarino
emphasized that an important ingredient for the IMEX R-K schemes 
of type A is $b_i = \tilde{b}_i$ for all $i$. Such a
choice provides a significant benefit for  the differential
component $u$, i.e., an order reduction does not appear for this
component. On the other hand, conditions \[e^T_s \tilde{A} =
\tilde{b}^T,\quad e^T_s A = b^T\] imply $a_{ss} = b_s \neq
\tilde{b}_s = \tilde{a}_{ss} =0$ which means that for a stiffly
accurate IMEX R-K scheme it is $b \neq \tilde{b}$, and therefore
we expect to observe order reduction for the differential variable.

\item It is impossible to construct a second
order stiffly accurate IMEX R-K scheme of type A with $s =
3$ internal stages. The proof is given in Appendix \ref{A2}. In practice,
in order to satisfy all these order conditions we have to increase
the number of the internal stages. In view of such difficulties,
for type A schemes, we shall consider second order IMEX R-K schemes
with $s = 3$ and $\tilde{b} = b$ in order to avoid the 
order reduction, giving
up to the FSAL property of the explicit scheme (and with that the global stiff accuracy of the IMEX scheme). 
An example is the scheme SSP(3,3,2) in Appendix \ref{A4}. In this case, if $q(u) \neq 0$, it is 
$g(u_{n+1},v_{n+1}) \neq 0$ as $\varepsilon \to 0$.

\item  Formulation (\ref{SPP}) in the limit case $\varepsilon \rightarrow
0$ yields the index-1 DAE. Then using the same technique adopted
in \cite{Bosc-09}, we can derive additional order conditions,
called \emph{algebraic conditions}, that guarantee the correct
behavior of the numerical solution in the limit $\varepsilon
\rightarrow 0$ and maintain the accuracy in time of the scheme.
If the implicit scheme is stiffly accurate, such conditions becomes, 
to various order of accuracy,
\begin{eqnarray}\label{oc_index1}
\begin{array}{l}
\tilde{c}_s  = 1, \ \textrm{(second order)} \\
e^T_{s}\tilde{A}\tilde{c}  = 1/2, \ \textrm{(third order)}
\end{array}
\end{eqnarray}
where $e_s = (0,...,0,1)^T$. If the IMEX schemes is globally stiffly accurate, then (\ref{oc_index1}) are automatically satisfied,
since $e^T_s \tilde{A} = \tilde{b}^T$.
\item Finally we observe that, in order to construct an order $p \geq 3$ IMEX R-K of type A 
and to maintain accuracy we have to increase the number of the classical order conditions too. 
Usually several simplifying assumptions (see \cite{Bosc-09}, \cite{Bosc-10},
\cite{HW} for details) could help to reduce the number of such
conditions, but, higher orders type A schemes are more
complicated to construct than CK or ARS schemes because of
additional order conditions (see \cite{Bosc-10}) due to the fact that $c\neq \tc$.
\end{itemize}

\subsection{Analysis of TYPE CK schemes}
Similar considerations about BPR approach, explained for the IMEX R-K scheme of type A in the limit case
$\varepsilon \to 0$, can be reproposed here
for the type CK when applied to the system (\ref{DA1}), with slightly modifications.
Of course, if we consider the general system (\ref{I3b}) we obtain again 
an IMEX R-K scheme of type CK in the diffusion limit, i.e. $\varepsilon \to 0$, where the diffusion term $f_2(u)$ is treated implicitly
and a CFL hyperbolic condition for the time step is required. 

Indeed, we consider an IMEX R-K schemes of type CK where, by Definition \ref{Def:CK},
we assume that the submatrix $\hat{A}$ is invertible and $a_{11} =
0$. The Butcher tableaux of a CK scheme takes the form
\begin{displaymath}
\begin{array}{c|c}
0         &0 \ \ 0\\
\hat{c} &a \ \ \hat{A}\\
\hline
 & b_1 \ \ \hat{b^T} \end{array} \quad
\end{displaymath}
with $a  = (a_{21},\ldots,a_{s,1})^T$ and
$\hat{b}^T=(b_{2},\ldots,b_{s})$. In order to simplify the
analysis we consider that the implicit part of the scheme is
stiffly accurate. Under this circumstance it is easy to prove that

\begin{eqnarray}\label{prorCK}
b_1 + \hat{b}^T\alpha = 0,
\end{eqnarray} 
where $\alpha \equiv - \hat{A}^{-1}a$
(see \cite{Bosc-07} for details).

Then, considering a scheme of the type CK, the second equation in (\ref{stypeA}) becomes% -(\ref{typeA}) can be written in the form
\begin{eqnarray}\label{secondeq}
\varepsilon^2 V_{k}   =  \varepsilon^2 v_n +  \Delta t a_{k1}g(u_n,v_n) + \Delta t\sum_{j=2}^{k} a_{kj} g(U_j,V_j). 
\end{eqnarray}
with $k = 2,...,s$.

Now multiplying by $\hat{\omega}_{kj}$, where $\hat{\omega}_{kj}$
are the elements of the inverse of $\hat{A}$, and summing on $k$,
we obtain
\begin{eqnarray}
\Delta t g(U_k,V_k) = \varepsilon^2\sum_{j=2}^{s}
\hat{\omega}_{kj}(V_j-v_n) + \Delta t \alpha_k g(u_n,v_n), \ \
\textrm{for} \ \ k=2,\ldots,s \nonumber
\end{eqnarray}
where \[\sum_{l = 2}^{s} \hat{\omega}_{kl}a_{lj} = \delta_{kj},
\qquad -\sum_{l = 2}^{s}\hat{\omega}_{kj}a_{j1} = \alpha_k.\]
Inserting the expression $\Delta t g(U_k,V_k)$ into the second
equation in (\ref{typeA}) we obtain
\begin{eqnarray}\label{index1}
%  u_{n+1} & = &  u_n  +  \Delta t \left(\tilde{b}_1 f_1(u_n,v_n) +  \sum_{k = 2}^{s} \tilde{b}_k  f_1(U_k,V_k) +  b_1 f_2(u_n) + \sum_{k = 2}^{s} b_k  f_2(U_k)\right) \nonumber         \acapo
 %             &   &  \acapo
 %           &   &                    \nonumber \acapo
\varepsilon^2 v_{n+1}  & = & \varepsilon^2 \mathcal{R}(\infty)  v_n +
\varepsilon^2\Delta t \sum_{k=2}^{s} b_k \omega_{kj}  V_j + \Delta
t \left(b_1 + \sum_{k=2}^{s} b_k\alpha_k\right)g(u_n,v_n)
%\nonumber
\end{eqnarray}
%with
%\begin{eqnarray}\label{index_1}
%   U_k & = &  u_n + \Delta t\left(\tilde{a}_{k1} f_1(u_n, v_n) + \sum_{j = 2}^{k-1} \tilde{a}_{kj} f_1(U_j,V_j) + a_{k1} f_2(u_n) + \sum_{j = 1}^{k} a_{kj}
%   f_2(U_j)\right).
%\end{eqnarray}
Then by (\ref{prorCK}) the last term in the second equation in (\ref{index1}) drops and in the limit case for $\varepsilon = 0$ we can write
\begin{eqnarray}
   v_{n+1}  & = & \mathcal{R}(\infty)  v_n +  \Delta t \sum_{k=2}^{s} b_k \omega_{kj}  V_j,      \nonumber
\end{eqnarray}
with
\begin{eqnarray}\label{pro2CK}
  g(U_k, V_k) & = & \alpha_k g(u_n,v_n), \ \ \textrm{for} \ \ k = 2,\ldots,s.  
\end{eqnarray}
Note that, for IMEX R-K schemes of type CK, the stability function
$\mathcal{R}(z)$ of the implicit part of the scheme takes the form
\begin{eqnarray}\label{STAB}
\mathcal{R}(z) &=& 1 + z(b_1 + \hat{b}^T(I-z\hat{A})^{-1}(\one_{s-1} + za)) \nonumber \acapo
& & \acapo
&=& (b_1 -\hat{b}^T \hat{A}^{-1}a)z + (1-\hat{b}^T\hat{A}^{-1}\one_{s-1} + \hat{b}^T\hat{A}^{-2}a) + \mathcal{O}(\frac{1}{z}). \nonumber
\end{eqnarray}
We obtained this result, by applying one step of the implicit part
of the scheme to the test problem $y' = \lambda y$, $y(t_0) = 1$,
with $\lambda \in \mathbb{C}$ and $\one_s = (1,\ldots,1)^T \in
\mathbb{R}^s$.

Thus, the only stiffly accurate condition, i.e. $\hat{e}_{s-1}^T \hat{A} = \hat{b}^T$ is not enough to guarantee that $lim_{z\rightarrow \infty}\mathcal{R}(z) = 0$ and then an additional condition is required for the implicit part of the scheme, (for details see \cite{Bosc-10}). This is expressed by the following
\begin{proposition} If
\begin{eqnarray}\label{eq:L}
-\hat{e}_{s-1}^{T}\hat{A}^{-1}a = \sum_{j\geq 2}\hat{\omega}_{sj}a_{j1} = 0,
\end{eqnarray}
then $\mathcal{R}(\infty) = 0$, where $\hat{e}_{s-1} = (0,\ldots,0,1)^{T}
\in \mathbb{R}^{s-1}$.
\end{proposition}\\
{\bf Proof.} In fact, assuming $\hat{A}$ invertible, we get $\hat{b}^T \hat{A}^{-1} = \hat{e}^T_{s-1}$ and when $z\rightarrow \infty$, from (\ref{STAB}) we obtain $\mathcal{R}(\infty) = \hat{b}^T \hat{A}^{-2}a = -\hat{e}_{s-1}^{T}\hat{A}^{-1}a$, which is zeros if (\ref{eq:L}) is satisfied. $\Box$ \\
Note that the previous Lemma implies that $\alpha_s = -\hat{e}^T_{s-1} \hat A^{-1} a = 0$ and by (\ref{pro2CK}) with $k = s$ we obtain $g(U_s,V_s)=0$, then the last stage values lie on the manifold  $g(u,v) = 0$ as $\varepsilon \rightarrow 0$. Now we observe that if the IMEX R-K scheme of type CK is globally stiffly accurate, we obtain from (\ref{index1}) and (\ref{index1A}) $u_{n+1} = U_s$ and $v_{n+1} = V_s$ and therefore $g(u_{n+1},v_{n+1}) = 0$ with $v_{n+1} = G(u_{n+1})$.

Since an IMEX R-K schemes of type ARS is a particular case of the
type CK where the vector $a = 0$, then the same results hold
true.

\paragraph{General remarks for type CK}
\begin{itemize}
\item IMEX CK schemes \cite{CK} are attractive because of their good
properties.  The implicit part of this scheme is \emph{singly
  diagonally implicit Runge-Kutta} (SDIRK) with $a_{ii} = \gamma>0$
for $i = 2,...,s$ and differs from the classical SDIRK one because
$a_{11} = 0$. In \cite{CK} such implicit schemes are called \emph{explicit}
singly diagonally implicit (ESDIRK). A consequence to set $a_{11} = 0$
is the possibility to guarantee stage-order $q$
higher than the in the case of SDIRK, for which $q = 1$.
Moreover here we consider schemes that are stiffly accurate according
to Definition \ref{1def}. Such schemes will project the solution on
the manyfold in the limit of infinite stiffness. 
For these schemes $b \neq \tilde{b}$, so one of the so-called simplifying
conditions cannot be applied \cite{Bosc-10}. 
Here we require that $c_{i} = \tilde{c}_i$ for all $i= 2,...,s$; this
choice will reduce the number of coupled order conditions. 
\end{itemize}

\section{IMEX-Finite Difference schemes} \label{Sec4}
When constructing numerical schemes, one has also to take a great
care in order to avoid spurious numerical oscillations arising
near discontinuities of the solution. This is avoided by a
suitable choice of space discretization. To this aim it is necessary to use 
non-oscillatory interpolating algorithms,
in order to prevent the onset of 
spurious oscillations (like ENO and WENO methods), see \cite{Shu}.
Moreover the choice of the space discretization may be relevant 
for a correct treatment of the boundary conditions.

In this section we emphasize some requirements about the space
discretization of the system (\ref{I3b}). We remark that the
dissipative nature of upwind schemes \cite{NP, NP2} depends
essentially on the fact that the characteristic speeds of the
hyperbolic part are proportional to $1/\varepsilon$. On the other
hand central differences schemes avoid excessive dissipation but
when $\varepsilon$ is not small or when the limiting equations
contain advection terms may lead to unstable discretizations. In
order to overcome these well-known facts and to have the correct
asymptotic behavior we fix some general requirements for the
space discretization.
\begin{enumerate}
 \item{\em Correct diffusion limit.} Let us consider system (\ref{I3b}) with $q(u) = 0$.
 In the limit case $\varepsilon \rightarrow 0$ the therm $v + \partial_x p(u) \rightarrow 0$ from the second equation.
 If we want that $v + \mu(\varepsilon) \partial_x p(u) \rightarrow 0$ also in the first equation, we need to use the same space discretization for
 the term $\partial_x p(u)$ and require that $\mu(0) = 1$.
    \item{\em Compact stencil.} Among the advantages of our approach there is the possibility to have a scheme with a compact stencil in the diffusion limit $\varepsilon \rightarrow 0$. This property is satisfied if point 1)
    is satisfied and we use a suitable discretization for the second order derivative that characterize the diffusion limit.
    \item{\em Shock capturing.} The schemes when $q(u)\neq 0$ should be based on shock capturing high order fluxes for the convection part. This is necessary not only for large values of $\varepsilon$ but also when we consider convection-diffusion type limit equations with small diffusion. The high
    order fluxes are then necessary for all space derivatives except for the second order term $\mu(\varepsilon) \partial_{xx} p(u)$ on the right-hand side.
\item{\em Avoid solving nonlinear algebraic equations.} In order to achieve this the implicit space derivative $\partial_x p(u)$ in the second equation
must be evaluated using only nodal values of $u$ which can be obtained from the solution of the first equation.
\end{enumerate}
The above properties are satisfied for example using high accuracy
in space obtained by finite difference discretization with
Weighted-Essentially Non Oscillatory (WENO) reconstruction,
\cite{Shu}. 

System (\ref{I3b}) may be written in the form
\begin{eqnarray}\label{hyper}
u_t +(v+\mu p(u)_x)_x &=& \mu p(u)_{xx},\nonumber\\
 & & \\
v_t  &=& \frac{1}{\varepsilon^2}\left(q(u) - (v + p(u)_x)\right). \nonumber
\end{eqnarray}
with $\mu = \mu(\varepsilon)$ introduced in Section 2.
The terms on the right-hand side will be treated implicitly.
For large value of $\varepsilon$ the explicit flux is just $(v,0)^T$, 
while for small values of $\varepsilon$ it is $(v+p(u)_x,0)^T$.
Here we describe a finite difference WENO scheme for a system of the form 
\begin{eqnarray*}
U_t + F(U)_x = R(U),
\end{eqnarray*}
and apply it to the system (\ref{hyper}) with 
\begin{eqnarray*}
F(U) = (v+ \mu p(u)_x,0)^T,\\
R(U) = (\mu p(u)_{xx},\frac{1}{\varepsilon^2}\left(q(u) - (v + p(u)_x)\right)).
\end{eqnarray*}
As $\varepsilon \rightarrow \infty$ and $\mu \rightarrow 0$, the system becomes
\begin{eqnarray*}
u_t + v_x &=& 0,\\
v_t       &=& 0
\end{eqnarray*}
and the characteristic speed of the system is $\lambda = 0$ (twice). As $\varepsilon \rightarrow 0$ and $\mu \rightarrow 1$,
$v + \mu p(u)_x \rightarrow q(u)$ and the system relaxes to the equation
\begin{eqnarray*}
u_t + q(u)_{x} = p(u)_{xx}
\end{eqnarray*}
and the characteristic speed of the left hand side is given by $\lambda = q'(u)$.

Conservative finite difference for system (\ref{hyper}) are written as follows, \cite{A}
\begin{eqnarray*}
\frac{d U_j}{dt} &=& -\frac{\hat{F}_{j+\frac{1}{2}}-\hat{F}_{j-\frac{1}{2}}}{\Delta x} + G(U_j)
\end{eqnarray*}
where $U_j(t) \approx U(x_j, t)$ is an approximation of the pointwise value of $U$ at grid nodes, 
and the numerical flux at cell edge $x_ {j+\frac{1}{2}}$ is computed as follows
\begin{eqnarray*}
\hat{F}_{j+\frac{1}{2}} = \hat{F}^{+}_{j}(x_{j+\frac{1}{2}})+\hat{F}^{-}_{j+1}(x_{j+\frac{1}{2}}).
\end{eqnarray*}
The function $\hat{F}_j^+(x)$ and $\hat{F}_{j+1}^-(x)$ are suitable reconstructions defined, respectively,
in cell $j$ and in cell $j+1$. They are obtained as follows. First, we assume that the flux can be split
into a positive and negative component 
\begin{eqnarray*}
F(U) = F^+(U) + F^-(U),
\end{eqnarray*}
with $\lambda(\nabla_U F^+(U)) \geq 0$, $\lambda(\nabla_U F^-(U)) \leq 0$.
The quantity $F^{\pm}_j = F^{\pm}(U_j)$ are computed at cell center.
Then $\hat{F}^{\pm}_j(x)$ are reconstructed from $\left\{F^{\pm}_j\right\}$
using high order essentially non oscillatory reconstruction, such as ENO or 
WENO, that allows pointwise reconstruction of a function from its cell averages, 
(see, e.g. \cite{Shu} for details).

The flux $F$ may contain derivatives. For example the first equation
in system (\ref{hyper}) contains $p(u)_x$. Such terms are computed by
point-wise WENO reconstruction.

In all our examples we used the simple local Lax-Friedrix flux decomposition, i.e.
$F^+(U) = \frac{1}{2}(F(U)+\alpha U)$, $F^-(U) = \frac{1}{2}(F(U)-\alpha U)$,
$\alpha \geq \textrm{max}_U\left|\rho(\nabla_U F)\right|$, $\forall A \in \mathbb{R}^{m \times m}$, where $\rho(A) = \textrm{max}_{1 \leq i \leq m} \left|\lambda_i(A)\right|$
denotes the spectrum radius of matrix $A$, and the max defining $\alpha$ is taken for $U$ varying in a suitable
range in a neighborhood of each cell.
In our test case we chose $\alpha = 1$, since as $\varepsilon \rightarrow \infty$, $\rho(\nabla_U F) = 0$ and 
in our numerical test $q(u)$ is either $0, \ u,$ or $u^2/2$, with $U$ ranging in $[0,1]$, therefore $\left|q'(u)\right|\leq1$.

We remark here that the choice of $\alpha$ is based on the
characteristic speeds of the limit convection-diffusion equation,
while a more detailed analysis is needed to justify its use in
intermediate regions, for which the characteristic speeds can be much
higher, and the stabilization that compensates for the apparent
violation of the hyperbolic CFL condition comes from the implicit
treatment of the diffusion term.

Furthermore for large value of $\varepsilon$, (e.g.,
$\varepsilon=1$), we want to avoid adding and subtracting terms
which may cause loss of accuracy.  For a semidiscrete scheme the function $\mu$ will 
depend also on the grid space $\Delta x$. A simple choice for $\mu$ is given by
\begin{eqnarray*}
  \mu = \exp(-\varepsilon^2/\Delta x)
 %   \mu(\varepsilon) =\begin{cases}  1,  & \mbox{if} \ \ \ \varepsilon < \Delta x,\\
  %                               0,  & \mbox{if}\ \ \ \varepsilon \geq \Delta
  %                               x,
%\end{cases}
\end{eqnarray*}
which is what we used in all our numerical tests.

For the diffusion term
$p(u)_{xx}$ we used the standard 2-nd order finite
difference technique for second order time discretization, and
the standard 4-th order finite difference technique where 3-rd
order time discretization are used.

\section{Numerical examples}\label{Sec5}
In this section we test several second and third order IMEX R-K
schemes presented in the literature that satisfy the algebraic order conditions (\ref{oc_index1}) and conditions in Definition 3.
Usually, IMEX time discretization are identified by an acronym (e.g. the
initials of the authors), and three numbers ($\sigma_E$, $\sigma_I$, $p$) denoting, respectively, the effective number of
stages (in practice the number of function evaluations) of the explicit and implicit scheme and the classical order of accuracy.

Below we list the IMEX R-K schemes used in the numerical tests.

\begin{itemize}
\item SSP(3,3,2): derived by Pareschi, Russo \cite{A},  it is a second order IMER R-K of type A, the explicit part is strongly stability preserving, while the implicit part is stiffly accurate.  In accordance with the proposition \ref{prop} in the Appendix \ref{A2}, this scheme is not globally stiffly accurate according to Definition \ref{1def}.

\item ARS(2,2,2): derived by Asher, Ruuth, Spiteri \cite{ARS}, it is a second order scheme, 
the double Butcher \emph{tableau} of this scheme is reproduced in Appendix \ref{A4}.
%This type of schemes are characterized by a triplet  $(s,\sigma,p)$  where $s$ is the number of stages of the implicit scheme,
%$\sigma$ is the number of explicit scheme stages with $\sigma = s+1$  if $\tilde{b}_i = b_i$ for all $i$,
%or $\sigma=s$ if $\tilde{b}_{s+1} = 0$. 
Note that this scheme is globally stiffly accurate according to the Definition \ref{1def}  
and satisfies the additional order conditions (\ref{oc_index1}).

\item ARS(4,4,3): derived by Asher, Ruuth, Spiteri \cite{ARS}, it is a third order scheme, 
the double Butcher \emph{tableau} of this scheme is reproduced in Appendix \ref{A4}. 
Similarly to ARS(2,2,2) this scheme is globally stiffly accurate according to Definition \ref{1def} 
and satisfies the additional order conditions (\ref{oc_index1}).

\item BPR(3,5,3) introduced in this paper is a third order IMEX R-K scheme of type CK and globally stiffly accurate according to Definition \ref{1def}. %i.e.  $\tilde{b}^T = e^T_s\tilde{A}$ and $b^T = e^T_s A$.
This scheme has $s = 5$ internal stages, $\sigma_E = 3$ explicit stages and $\sigma_I = 5$ implicit stages. The additional order conditions
(\ref{oc_index1}) are satisfied. This scheme is more efficient than ARS(4,4,3) for the explicit part, but less efficient for the implicit one. In many cases the computation of the explicit term is more expensive than the solution of the implicit step, thus resulting in an overall improvement in efficiency per time step. Furthermore, the coefficients of the explicit scheme are all non negative, which is an advantage for the stability properties of the scheme.
We reproduced the coefficients of this scheme in Appendix \ref{A4}. 
\end{itemize} 
%%%%%%%%%%%%%%%%%%%%%%%%%%%%%%%%
In all the computations presented in this paper we denote each scheme
with an acronym indicating the IMEX scheme and the type of space
discretization. 

Space discretization is identified by a short name containing the order of accuracy in space; for example, 
WENO53 (or WENO32), see for details \cite{Shu}, means a fifth (or a third) order reconstruction which reduced to third (or second) order near singularities and CdS2 stands for second order central discretization scheme. 

We remark here that all the analysis performed in the paper and the
numerical tests are performed under the assumption that the initial
data is {\em well-prepared}, which means that the initial condition 
lies in the limit manyfold as $\varepsilon\to 0$. If this condition is
not satisfied, then a loss of accuracy is observed, unless some
initial layer fix is adopted. Schemes of type A are more robust
against this problem, as is described in \cite{A}.

\subsection{Convergence test}
In this section we investigate numerically the convergence rate of the
second and third IMEX R-K schemes introduced before for a wide range
of the parameter $\varepsilon$. To this aim we apply these schemes to
simple prototype hyperbolic system (\ref{I3b}), with initial
conditions chosen in such a way that the exact solutions is smooth and
does not present a rapidly varying transient for small values of
$\varepsilon$. This is achieved in practice by imposing that the
initial condition satisfies the limit relation between $u$ and $v$ as 
$\varepsilon \to 0$. 

Numerical convergence rate is calculated by the formula
\begin{eqnarray*}
p = \log_2(E_{\Delta t_1}/E_{\Delta t_2}),
\end{eqnarray*}
where $E_{\Delta t_1}$ and $E_{\Delta t_2}$ are the global errors
computed with step $\Delta t_1 = \mathcal{O}(\Delta x)$,  and 
$\Delta t_2 = \Delta t_1/2$.
In the following tests we put $\varepsilon^2 = 10^{-6}$ and we choose
$\Delta t \approx \Delta x$.

For the first test we set $p(u) = u$
and $q(u) = 0$. Then we get
\begin{eqnarray}\label{lde}
\begin{array}{ll}
u_t = -v_x - \mu u_{xx} + \mu u_{xx},\\
\varepsilon^2 v_t = -u_x - v,
\end{array}
\end{eqnarray}
that in the limit case, $\varepsilon = 0$ and $\mu = 1$ leads to the linear diffusive problem
\begin{eqnarray}\label{Limdiff}
u_t = u_{xx},  \ \ \ u(x,0) = u_0(x).
\end{eqnarray}
We use periodic boundary conditions with $u_0(x) = \cos(x)$, and
$x \in [0 \ , 2\pi]$, so that $u(x,t) = u_0(x)\exp(-t)$ is an exact
solution of (\ref{Limdiff}). The final time is $T = 1$ and $\Delta t = 0.5 \Delta x$.

The results are reported in Table 5.1 and 5.2 showing that the expected convergence rates are reached for the
$u$-component.

Next we set $p(u) = q(u) = u$ and consider the following system
\begin{eqnarray*}
  u_t + v_x & = & \mu u_{xx} - \mu u_{xx}                  \nonumber        \acapo
          &  &                    \nonumber \acapo
 \varepsilon^2 v_t + u_x & = & - (v-u),      \nonumber
\end{eqnarray*}
where the limiting behavior is given by an advection-diffusion
equation. We use periodic boundary conditions with the initial
data $u(x,0) = \exp(-(1+\cos(x-\pi))/\sigma)$, $v(x,0) = u(x,0)(1 -
\mu\sin(x-\pi)/\sigma)$ with $\sigma = 0.05$ and $\mu =1$, on the
spatial interval $[0, 2\pi]$, at the final time $T = 0.3$ anfd $\Delta t = 0.5 \Delta x$. 
As reference solution we use the truncated Fourier representation of
the exact solution
\begin{eqnarray*}
U_{exa}(x,t) = \sum_{k = -\infty}^{+\infty} U_k(t) e^{ikx}, \ \ V_{exa}(x,t) = \sum_{k = -\infty}^{+\infty} V_k(t) e^{ikx}
\end{eqnarray*}
with $U_k(t)$ and $V_k(t)$ satisfying
\begin{eqnarray}\label{FC1}
\begin{array}{l}
\dot{U}_k = -i k V_k,\\
\varepsilon^2\dot{V}_k = -i k U_k +U_k -V_k.
\end{array}
\end{eqnarray}

For each $k$, system (\ref{FC1}) can be written as a $2 \times 2$ constant coefficient homogeneous system which can be solved exactly.
The results are given in Table 5.3 showing that
again the expected convergence rates are reached for the
$u$-component by all schemes.

\begin{table}%\label{tab1}
\caption{Convergence rate for $u$ in $L_\infty$-norm.}

\centering \footnotesize
\begin{tabular}{|c|c|c|c|}%\toprule
\hline
\lower.3ex\hbox{ARS(2,2,2)-CdS2} & \lower.3ex\hbox{SSP2(3,3,2)-CdS2} &\lower.3ex\hbox{  ARS(2,2,2)-WENO32}&\lower.3ex\hbox{SSP2(3,3,2)-WENO32} \\

\lower.3ex\hbox{$e^T_s \tilde{A} = \tilde{b}^T$}& \lower.3ex\hbox{$e^T_s \tilde{A} \neq \tilde{b}^T$}&\lower.3ex\hbox{$e^T_s \tilde{A} = \tilde{b}^T$}& \lower.3ex\hbox{$e^T_s \tilde{A} \neq \tilde{b}^T$}\\
 \lower.3ex\hbox{ $N$}  \  \lower.3ex\hbox{$L_{\infty}(u)$}  \ \lower.3ex\hbox{Order}& \lower.3ex\hbox{$N$}  \ \lower.3ex\hbox{$ L_{\infty}(u)$}  \ \lower.3ex\hbox{Order}&  \lower.3ex\hbox{$N$} \ \lower.3ex\hbox{$ L_{\infty}(u)$}  \ \lower.3ex\hbox{Order}  & \lower.3ex\hbox{$N$} \  \lower.3ex\hbox{$L_{\infty}(u)$} \  \lower.3ex\hbox{Order}   \\
\hline
\hline
   \lower.3ex\hbox{20}   \lower.3ex\hbox{7.800e-03} \ \ \ \ \ \ &\lower.3ex\hbox{20}  \lower.3ex\hbox{2.906e-02}\ \ \ \ \ \ &
   \lower.3ex\hbox{20}  \lower.3ex\hbox{4.820e-03} \ \ \ \ \ \  &  \lower.3ex\hbox{20}   \lower.3ex\hbox{4.697e-03} \ \ \ \ \ \ \\
                           
   \lower.3ex\hbox{40}  \lower.3ex\hbox{1.873e-04}  \lower.3ex\hbox{$2.05$}&
   \lower.3ex\hbox{40}   \lower.3ex\hbox{7.979e-03} \lower.3ex\hbox{$1.86$} &              
   \lower.3ex\hbox{40}   \lower.3ex\hbox{1.492e-03}  \lower.3ex\hbox{$1.69$} &      
   \lower.3ex\hbox{40}  \lower.3ex\hbox{1.483e-03} \lower.3ex\hbox{$1.66$}\\             
  
    \lower.3ex\hbox{80}   \lower.3ex\hbox{ 4.597e-04}  \lower.3ex\hbox{$ 2.02$}&
    \lower.3ex\hbox{80}   \lower.3ex\hbox{2.039e-03}   \lower.3ex\hbox{$1.96$}&
    \lower.3ex\hbox{80}   \lower.3ex\hbox{4.124e-04}  \lower.3ex\hbox{$1.85$}&
     \lower.3ex\hbox{80}   \lower.3ex\hbox{4.102e-04}  \lower.3ex\hbox{$1.85$}\\
     
      \lower.3ex\hbox{160} \lower.3ex\hbox{ 1.138e-04}  \lower.3ex\hbox{$2.01$}&
       \lower.3ex\hbox{160}  \lower.3ex\hbox{5.120e-04}   \lower.3ex\hbox{$1.99$}&
       \lower.3ex\hbox{160}  \lower.3ex\hbox{1.082e-04}   \lower.3ex\hbox{$1.93$}&
        \lower.3ex\hbox{160}  \lower.3ex\hbox{1.074e-04} \lower.3ex\hbox{$1.93$}\\
        
        \lower.3ex\hbox{320}  \lower.3ex\hbox{ 2.833e-05}   \lower.3ex\hbox{$2.00$}&
        \lower.3ex\hbox{320}   \lower.3ex\hbox{1.274e-04} \lower.3ex\hbox{$2.00$}&
         \lower.3ex\hbox{320}  \lower.3ex\hbox{2.760e-05} \lower.3ex\hbox{$1.97$}&
         \lower.3ex\hbox{320} \lower.3ex\hbox{2.748e-05} \lower.3ex\hbox{$1.96$}\\

                            \hline
\end{tabular}
\end{table}
%%%%%%%%
\begin{table}%\label{tab2}
\caption{Convergence rate for $u$ in $L_\infty$-norm.}

\centering \footnotesize
\begin{tabular}{|c|c|}%\toprule
\hline
\lower.3ex\hbox{ARS(4,4,3)-WENO53} &\lower.3ex\hbox{BPR(5,5,3)-WENO53}\\
\lower.3ex\hbox{$e^T_s \tilde{A} = \tilde{b}^T$}& \lower.3ex\hbox{$e^T_s \tilde{A}= \tilde{b}^T$}\\
 \lower.3ex\hbox{$N$}  \  \lower.3ex\hbox{$L_{\infty}(u)$}  \ \lower.3ex\hbox{Order}& \lower.3ex\hbox{$N$}  \  \lower.3ex\hbox{$ L_{\infty}(u)$}  \ \lower.3ex\hbox{Order}\\
\hline
\hline
 \lower.3ex\hbox{20}  \lower.3ex\hbox{1.810e-02} \ \ \ \ \ \ &  \lower.3ex\hbox{20}  \lower.3ex\hbox{1.639e-02} \ \ \ \ \ \   \\
 \lower.3ex\hbox{40}  \lower.3ex\hbox{3.365e-03} \lower.3ex\hbox{$2.42$} &      \lower.3ex\hbox{40}  \lower.3ex\hbox{3.099e-03}  \lower.3ex\hbox{$ 2.40$}\\
 \lower.3ex\hbox{80}  \lower.3ex\hbox{5.349e-04}  \lower.3ex\hbox{$2.65$}&  \lower.3ex\hbox{80}    \lower.3ex\hbox{5.167e-04}  \lower.3ex\hbox{$2.58$}\\
 \lower.3ex\hbox{160}  \lower.3ex\hbox{5.960e-05}  \lower.3ex\hbox{$3.16$}&\
           \lower.3ex\hbox{160}  \lower.3ex\hbox{5.821e-05} \lower.3ex\hbox{$3.14$}\\
 \lower.3ex\hbox{320}  \lower.3ex\hbox{5.968e-06}  \lower.3ex\hbox{ $3.31$}  & \lower.3ex\hbox{320}   \lower.3ex\hbox{5.949e-06}  \lower.3ex\hbox{$3.29$}\\
\hline
\end{tabular}
\end{table}

\begin{table}
\caption{Convergence rate for $u$ in $L_\infty$-norm in the
convection-diffusion limit.}

\centering \footnotesize
\begin{tabular}{|c|c|c|c|}%\toprule
\hline
\lower.3ex\hbox{ARS(2,2,2)-CdS2} & \lower.3ex\hbox{SSP2(3,3,2)-CdS2} & \lower.3ex\hbox{ARS(4,4,3)-WENO53} &\lower.3ex\hbox{BPR(5,5,3)-WENO53}    \\
\lower.3ex\hbox{$e^T_s \tilde{A} = \tilde{b}^T$}& \lower.3ex\hbox{$e^T_s \tilde{A} \neq \tilde{b}^T$}&\lower.3ex\hbox{$e^T_s \tilde{A} = \tilde{b}^T$}&\lower.3ex\hbox{$e^T_s \tilde{A} = \tilde{b}^T$}\\
\lower.3ex\hbox{$N$}  \ \lower.3ex\hbox{$ L_{\infty}(u)$}  \ \lower.3ex\hbox{Order}  & \lower.3ex\hbox{$N$} \  \lower.3ex\hbox{$L_{\infty}(u)$} \  \lower.3ex\hbox{Order}  & \lower.3ex\hbox{$N$} \ \lower.3ex\hbox{$L_{\infty}(u)$} \  \lower.3ex\hbox{Order} & \lower.3ex\hbox{$N$} \  \lower.3ex\hbox{$L_{\infty}(u)$} \  \lower.3ex\hbox{Order}  \\
\hline
\hline
   \lower.3ex\hbox{40}  \lower.3ex\hbox{3.867e-03} \ \ \ \ \ \ & \lower.3ex\hbox{40} \lower.3ex\hbox{2.615e-03} \ \ \ \ \ \ & \lower.3ex\hbox{40}  \lower.3ex\hbox{4.297e-04} \ \ \ \ \ \  &  \lower.3ex\hbox{40}  \lower.3ex\hbox{8.300e-04} \ \ \ \ \ \ \\
   
                             \lower.3ex\hbox{80}  \lower.3ex\hbox{9.457e-04}  $\lower.3ex\hbox{2.03}$ &
                              \lower.3ex\hbox{80}  \lower.3ex\hbox{6.243e-04}   \lower.3ex\hbox{$2.06$}&
                          \lower.3ex\hbox{80}  \lower.3ex\hbox{5.770e-05}  \lower.3ex\hbox{2.89}&
                            \lower.3ex\hbox{80}  \lower.3ex\hbox{1.167e-04}  \lower.3ex\hbox{2.83}\\
                          
                             \lower.3ex\hbox{160}  \lower.3ex\hbox{2.330e04}   \lower.3ex\hbox{2.02} &
                               \lower.3ex\hbox{160}   \lower.3ex\hbox{1.543e-04}  \lower.3ex\hbox{2.01}&
                             \lower.3ex\hbox{160} \lower.3ex\hbox{7.922e-06}   \lower.3ex\hbox{2.86}&
                                                 \lower.3ex\hbox{160} \lower.3ex\hbox{1.603e-05}  \lower.3ex\hbox{2.86}\\
                             
                              \lower.3ex\hbox{320}   \lower.3ex\hbox{5.798e-05}   \lower.3ex\hbox{2.00}&         
                            \lower.3ex\hbox{320}     \lower.3ex\hbox{3.850e-05}   \lower.3ex\hbox{2.00}&
                              \lower.3ex\hbox{320}  \lower.3ex\hbox{1.256e-06}   \lower.3ex\hbox{2.65} &
                             \lower.3ex\hbox{320}   \lower.3ex\hbox{2.230e-06}    \lower.3ex\hbox{2.85} \\
\hline
\end{tabular}
\end{table}
%%%%%%%%%%%%
The above convergence analysis has been performed in the limit
$\varepsilon \rightarrow 0$, therefore we might expect a
degradation of the accuracy for intermediate regimes as in the
case of  hyperbolic relaxation when the classical order is greater
then two \cite{Bosc-07, Bosc-10, CK}. Furthermore, from the
practical point of view, the understanding of this phenomenon is
essential in situations where one is interested in the
construction of higher order methods.
 
\begin{figure}[h]
\centering
\includegraphics[width=0.48\textwidth]{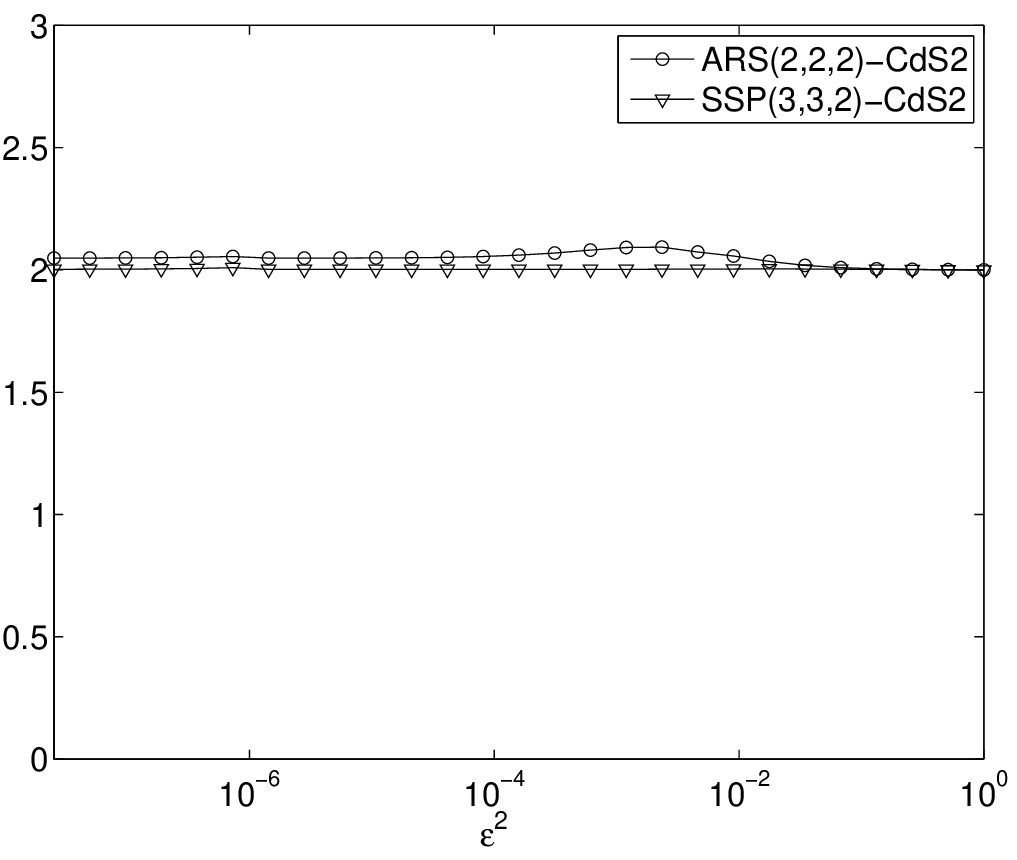}
\includegraphics[width=0.48\textwidth]{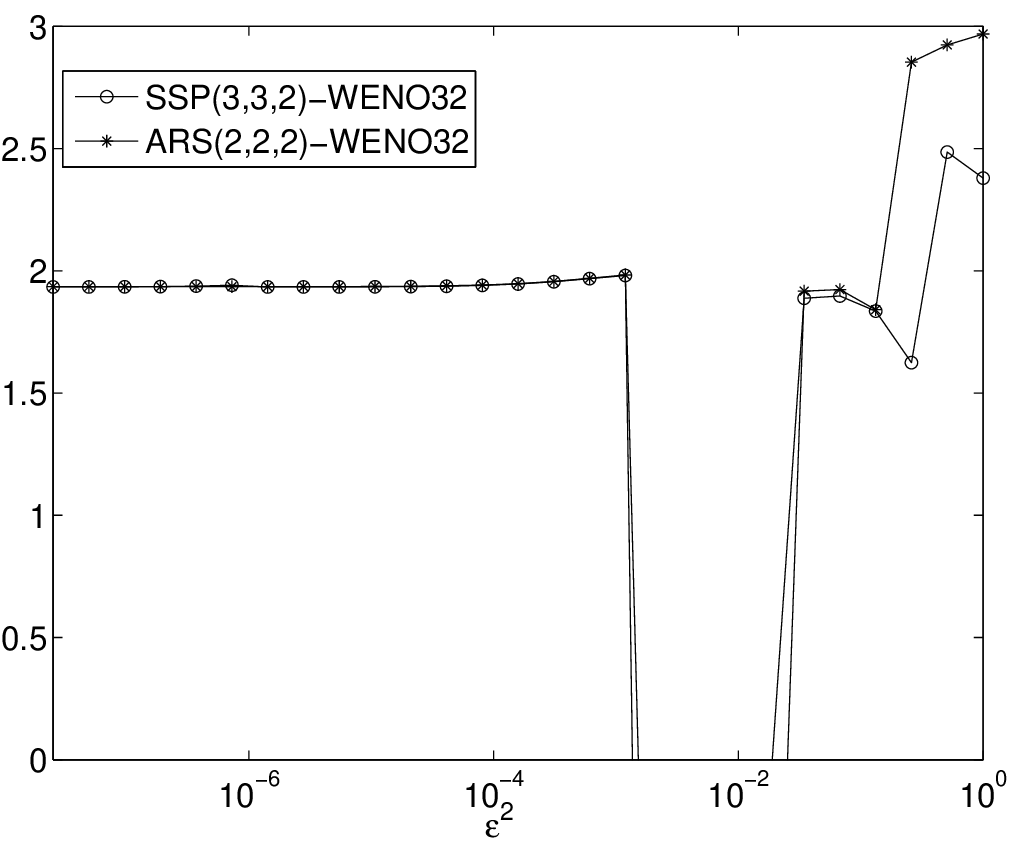}
\includegraphics[width=0.48\textwidth]{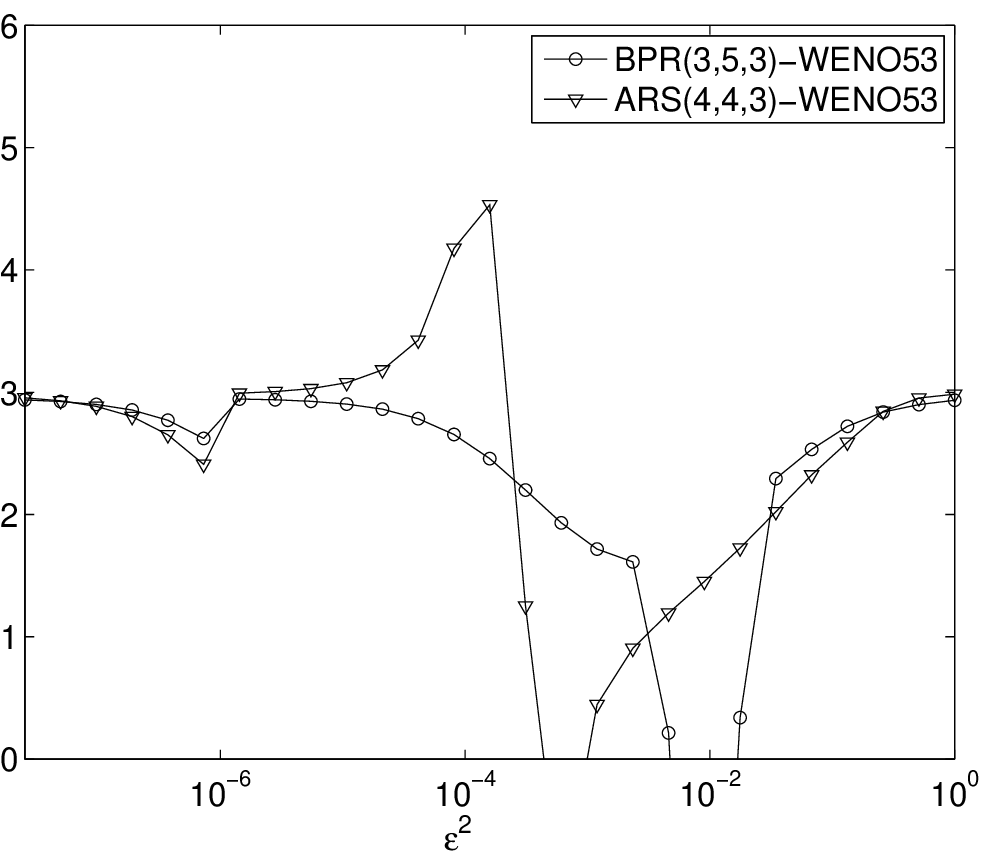}
\includegraphics[width=0.48\textwidth]{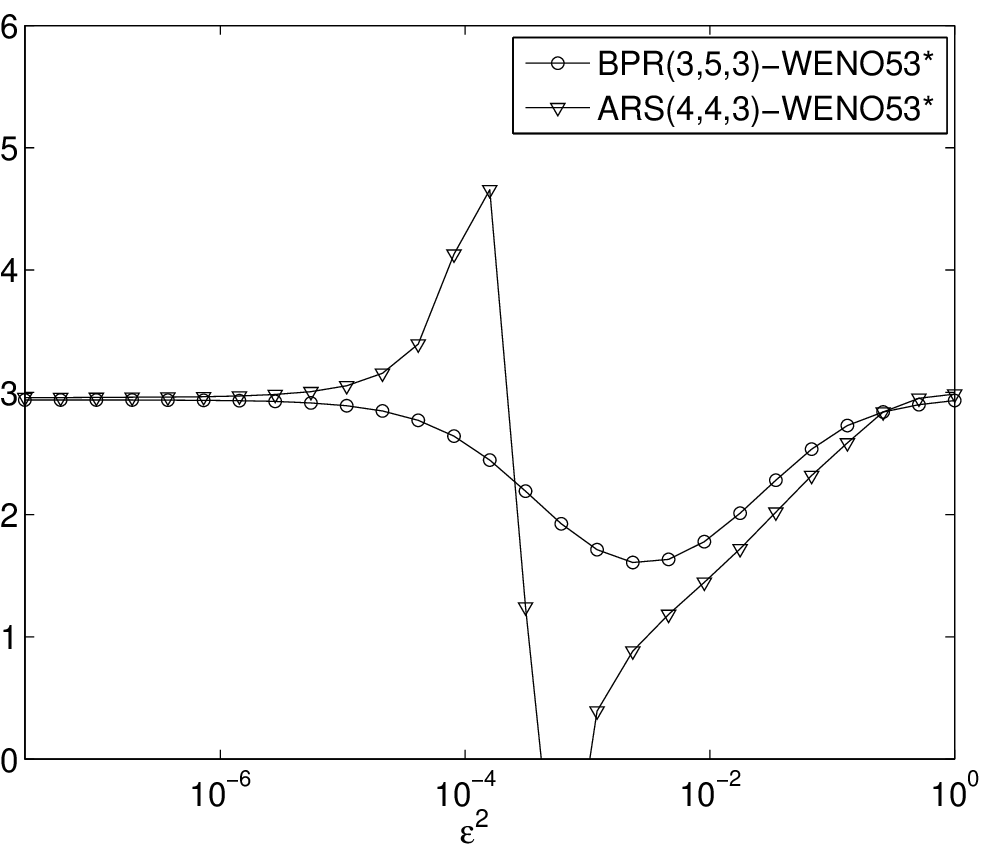}
\caption{Convergence rate in $L_\infty$-norm versus $\varepsilon^2$
  for various schemes. Lack of convergence for intermediate values of
  $\varepsilon$ is evident in the upper
  right and lower left panels. The lower right panel shows results
  obtained by a scheme in which the explicit and implicit term $\mu
  p(u)_{xx}$ is discretized spatially in identical way, namely by a
  second order discrete Laplacian.}
\label{fig:conv1}
\end{figure}

%Figures \ref{fig:conv1}, \ref{fig:conv2} and \ref{fig:conv3} show the
Figure \ref{fig:conv1} shows the
convergence rates as a function of $\varepsilon^2$ using different
values of $\varepsilon^2$ ranging from $10^{-6}$ to $1$ and $\Delta t \simeq \Delta x$. 
Second order schemes ARS(2,2,2)-CdS2, SSP2(3,3,2)-CdS2 have the prescribed order of accuracy
uniformly in $\varepsilon^2$ (upper left panel). 
%(see Fig. \ref{fig:conv1}). 
Instead, ARS(2,2,2)-WENO32 and SSP2(3,3,2)-WENO32 present a
degradation of accuracy at intermediate regimes (upper right panel). 

A similar lack of convergence in
  intermediate regimes is observed for both the third order schemes
  ARS(4,4,3)-WENO53 and BPR(3,5,3)-WENO53 (lower left panel). This results have a very
  different nature than the accuracy degradation observed in IMEX
  schemes applied to hyperbolic systems with stiff relaxation
  \cite{Bosc-07}.
  A plausible reason here appears to be that in intermediate regimes the two terms which
  are added and subtracted in the equations, i.e. $\pm \mu p(u)_{xx}$,
  are discretized in two very different ways: one is computed inside
  the flux (ARS(4,4,3)-WENO53 and BPR(3,5,3)-WENO53), and the other one (ARS(4,4,3)-WENO53* and BPR(3,5,3)-WENO53*)
  is discretized by a discrete one
  dimensional Laplacian, therefore the two terms do not almost cancel
  each other. Although in certain regimes such problem
  could be solved by treating the term $p(u)_x$ out of the flux (see,
  for example, the result in the lower right panel of Figure
  \ref{fig:conv1} for BPR(3,5,3)) and
  discretizing both terms $\pm \mu p(u)_{xx}$ in the same way, this
  may compromise the cancelation of the quantity $q(u)-v-p(u)_x$ in the
  flux. A general understanding and a robust treatment of intermediate
  regimes is beyond the scope of the present paper and requires
  further investigation.

%%%%%%%%%%%%%%%%%%%%%%%%%%%%%%%%%%%%%%%%%%%%%%%%%%%%%%%

\subsection{Shock test cases}

In this section we apply the scheme to problems with discontinuous
initial data, that in the limit as $\varepsilon \to 0$ reduce to
convection-diffusion equation.  Notice that in the relaxed limit the
scheme becomes an IMEX scheme for the limit equation. This test is
used to check both the shock capturing properties of the scheme, and
its relaxation to an IMEX scheme for the limit equation.

In the rest of the section we will consider the third
order BPR(3,5,3)-WENO53 scheme. 

First we consider a purely diffusive linear problem. We solve a
Riemann problem,  in the rarefied and diffusive regime for system
\begin{eqnarray*}
  u_t + v_x & = &  \mu u_{xx} -  \mu u_{xx}                  \nonumber        \acapo
            &   &                 \\
 \varepsilon^2 v_t + u_x & = & - v.      \nonumber
\end{eqnarray*}
We take the following initial data
\begin{eqnarray*}
u_L &=& 2.0 \ \ \ v_L = 0, \ \ -1 < x < 0,\\
u_R &=& 1.0 \ \ \ v_R = 0, \ \ 0 < x < 1.
\end{eqnarray*}
As $\varepsilon$ goes to zero we get $u_t = u_{xx}$, i.e. the
problem becomes a classical Riemann problem for the heat equation.

In order to test our scheme we compute the numerical solution in
the rarefied ($\varepsilon^2 > \Delta x$) regime and in the
diffusive ($\varepsilon^2 \ll \Delta x$) regime. This means that
when $\varepsilon^2$ is very large (i.e., rarefied regime) $\mu$ is
very small, and on the other hand when $\varepsilon$ is very small
(i.e., diffusive regime) $\mu$ is equal to $1$.

\begin{figure}%[ht]
\centering
%\vspace{2.5in}
\includegraphics[width=0.45\textwidth]{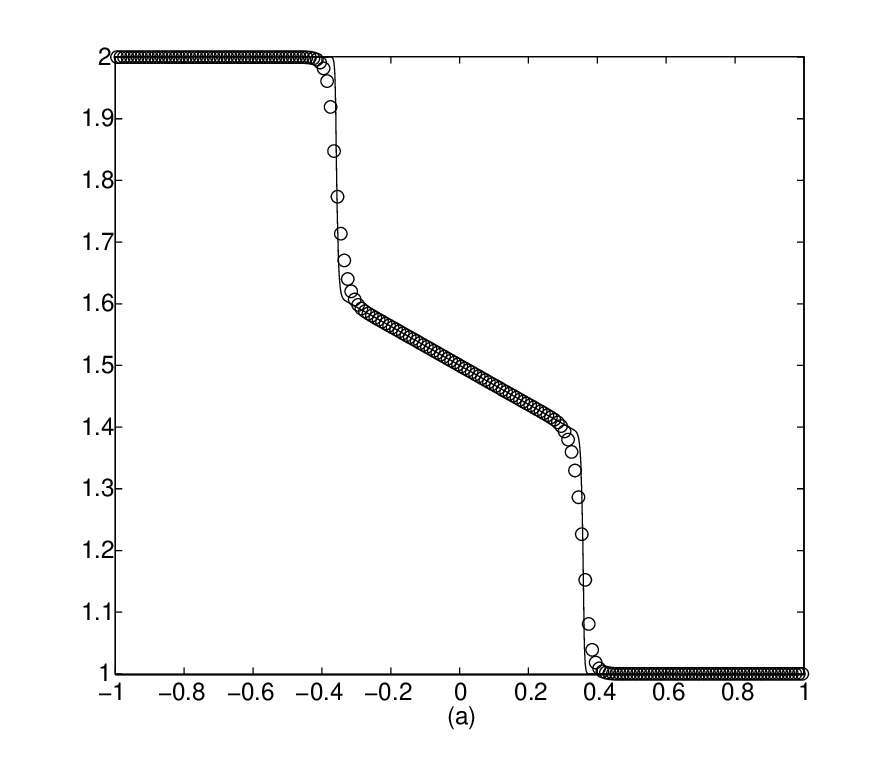}
\includegraphics[width=0.45\textwidth]{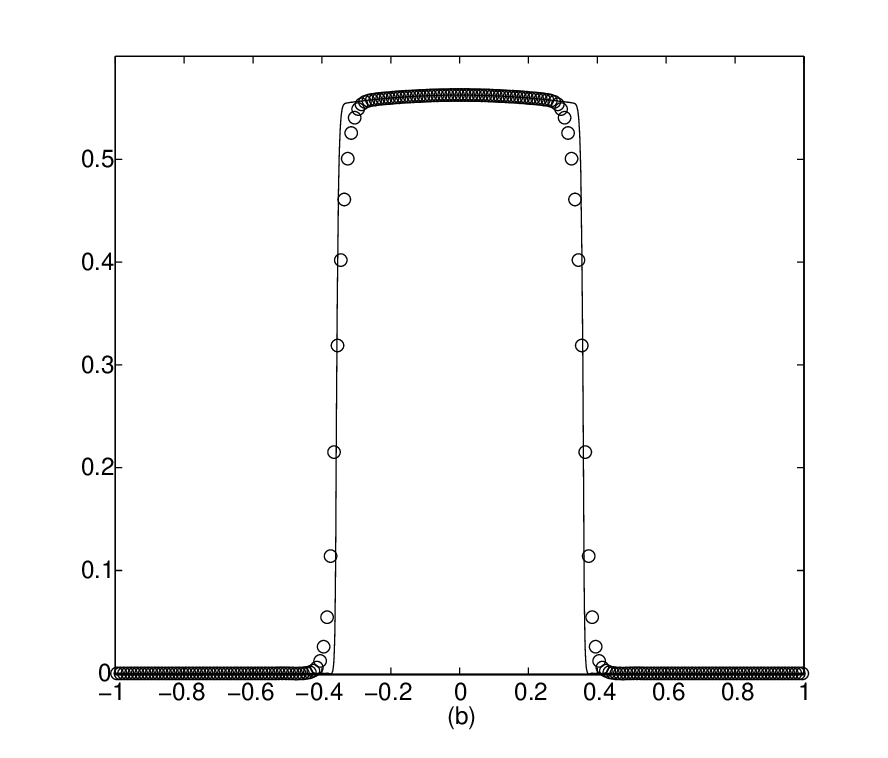}
\caption{Numerical solutions at time $t = 0.25$ in the rarefied
regime ($\varepsilon^2 = 0.7$) with $\Delta t= 0.5\Delta x$ and $\Delta x = 0.01$. On the left-hand side the
mass density $u$ (a) and on the right-hand side the flow $v$ (b). Solid line is the
reference solution.} \label{fig:1}
\end{figure}

\begin{figure}\label{limreg}%[ht]
\centering
%\vspace{2.5in}
%\includegraphics[width=0.8\textwidth]{RHO_LIMIT2.eps}
\includegraphics[width=0.45\textwidth]{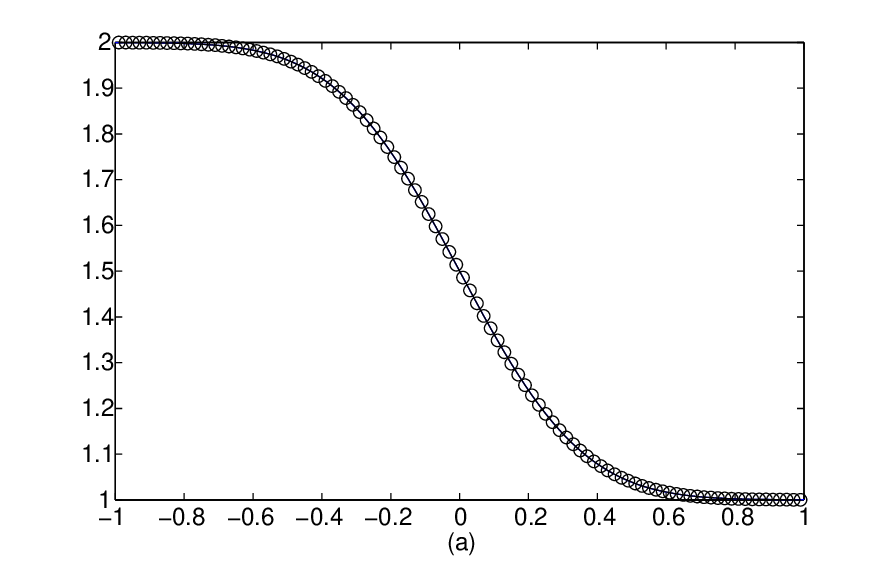}
\includegraphics[width=0.45\textwidth]{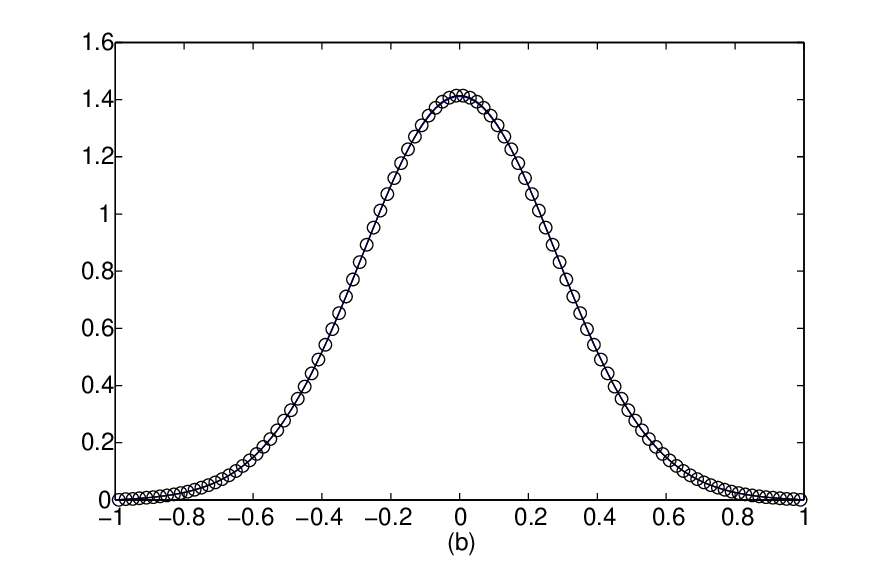}
\caption{Numerical solutions at time $t = 0.04$ in the parabolic
regime ($\varepsilon^2 = 10^{-6}$) with $\Delta t = 0.5\Delta
x$ and $\Delta x = 0.02$. On the left-hand side the mass density (a) $u$ and on the right-hand side the flow $v$ (b). Solid line
is the reference solution.} \label{fig:2}
\end{figure}

We set $\varepsilon^2 = 0.7$ for the rarefied regime  and $\varepsilon^2
= 10^{-6}$ for the diffusive regime (or stiff regime). 
The numerical solution for $u$ and $v$ in the rarefied (Fig.\ref{fig:1})
and diffusive regime (Fig.\ref{fig:2}) are depicted with a reference solution
obtained using a fine spatial grid of $N = 2000$ cells. 

As boundary conditions we set $v = 0$ at $x = \pm 1$. Compatibility with the system gives
$u_x = 0$ at $x = \pm 1$. Notice that the characteristic variables for this problem are $\xi_{\pm} = u \pm \varepsilon v$,
therefore condition $v = 0$ at $x = \pm 1$ is equivalent to impose $\xi_{+} = \xi_{-}$ at the boundary.
For such a reason we denote these boundary conditions as \emph{reflecting} boundary conditions.  

The solution is reported at final time $t = 0.25$ in the rarefied regime (Fig. \ref{fig:1}) and $t =
0.04$ in the diffusive regime (Fig. \ref{fig:2}).
%%%%%%%
In the figures we observe that the scheme captures well the
correct behavior of the solutions both in rarefied regime
where it provides an accurate description of the shock without
oscillations near the discontinuities, and in the diffusive regime
where the numerical solution matches accurately the reference
solution.
\begin{figure}%[ht]
\centering
%\vspace{2.5in}
%\includegraphics[width=0.8\textwidth]{RWrho.eps}
%\includegraphics[width=0.8\textwidth]{RWrho_limit.eps}
\includegraphics[width=0.45\textwidth]{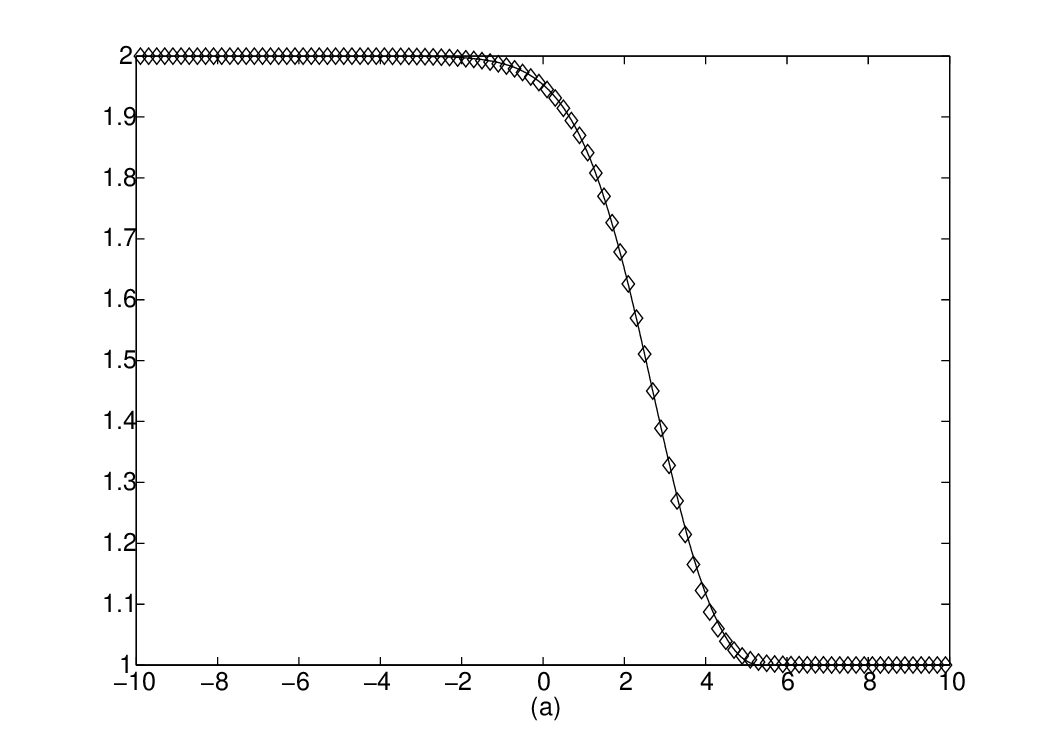}
\includegraphics[width=0.45\textwidth]{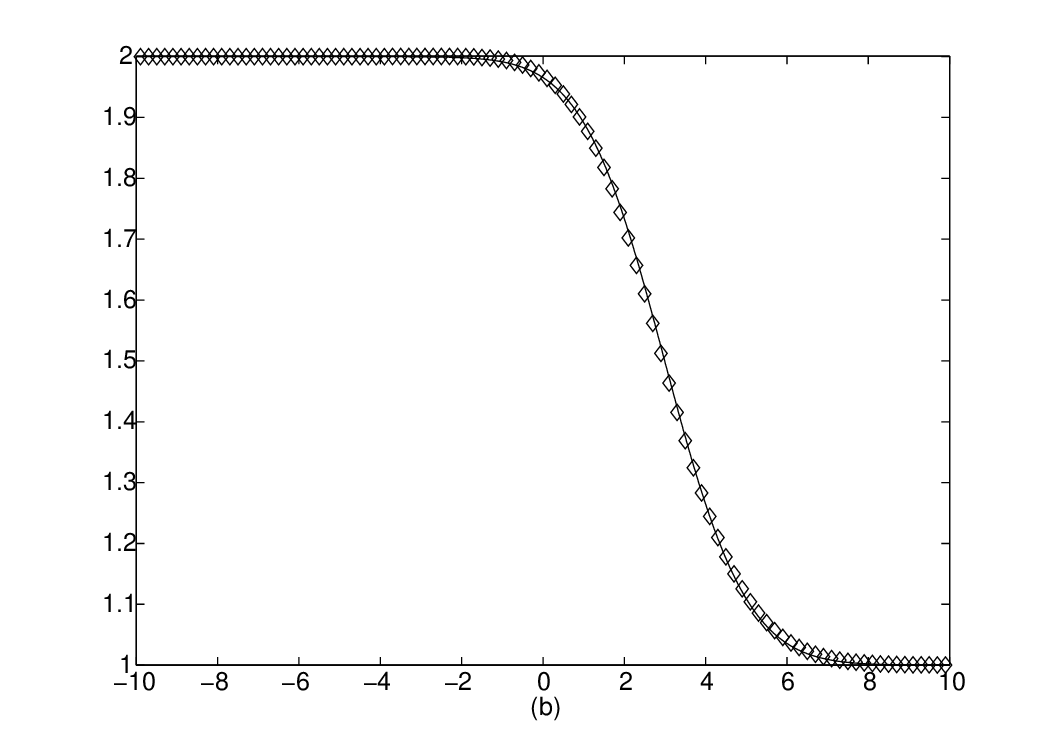}
\caption{Numerical solutions at time $t = 2.0$, (a) in the intermediate regime ($\varepsilon^2 = 0.4$) and (b) parabolic regime
($\varepsilon^2 = 10^{-6}$) with $\Delta t = 0.25\Delta x$ and $\Delta x = 0.2$. Solid line is the exact
solution.} \label{fig:3}
\end{figure}

Finally we consider the nonlinear Ruijgrok-Wu model, \cite{RW}, (for details see \cite{JPT})
\begin{eqnarray}\label{rw}
  u_t + v_x & =& \mu(\varepsilon) \frac{u_{xx}}{2 k_0} - \mu(\varepsilon) \frac{u_{xx}}{2 k_0}                    \nonumber        \acapo
            &   &                     \\
 \varepsilon^2 v_t + u_x & = & - 2k_0\left[v - \frac{C}{2}(u^2 - \varepsilon^2 v^2)\right]   \nonumber
\end{eqnarray}
where we added and subtracted to the original model the quantity $\mu u_{xx}$ on the right-hand side of the first equation
and $C$ is a constant, (in our test we chose $C = 1$). In the diffusive limit $\varepsilon
\rightarrow 0$, the second equation provides
\begin{eqnarray*}
v = \frac{1}{2}u^2 -\frac{1}{2 k_0}u_x,
\end{eqnarray*}
and we get the limiting viscous Burgers equation
\begin{eqnarray}\label{Burg}
  u_t +  \left(\frac{u^2}{2}\right)_x & =  \frac{u_{xx}}{2 k_0}.
\end{eqnarray}
%with $\nu = 1/(2 k_0)$. 
Note that through BPR approach
the IMEX R-K scheme in the diffusive limit relaxes to the same IMEX R-K one
for equation (\ref{Burg}) where the convection term is treated explicitly and the diffusion term implicitly. 

The exact solution to the shock-wave problem has
been given in \cite{RW}. The initial conditions are two local
Maxwellians characterized by
\begin{eqnarray*}
u_L = 2.0, \ \ \ -10<x<0,\\
u_R = 1.0, \ \ \ 0<x<10,
\end{eqnarray*}
with $v=\left[(1+u^2\varepsilon^2)^{1/2}-1\right]/\varepsilon^2$. 

In Figure \ref{fig:3} we show the computed solution for the mass
density $u$ in the intermediate ($\varepsilon^2 = 0.4$) and
parabolic ($\varepsilon^2 = 10^{-6}$) regimes versus the exact
solution. As it can be seen once again, the scheme gives an
accurate description of the viscous shock profiles.
%%%%%%%
%%%%%%%%%%%%%%%%%%%%%%%5

\section{Application to transport equations}\label{Sec6}
In this section we apply the IMEX schemes derived in the first
part of the manuscript to the case of neutron transport equations
\cite{Klar, JPT2, LM}.

We consider the multidimensional transport equation
under the diffusive scaling. Let $f(t,\textbf{x},\textbf{v})$ be
the probability density distribution for particles at space point
$\textbf{x} \in \mathbb{R}^d$, time $t$, traveling with velocity
$\textbf{v} \in \Omega \subset \mathbb{R}^d$ with $\int_{\Omega}
d\textbf{v} = S$. Here $\Omega$ is symmetric in $\textbf{v}$,
meaning that $\int_{\Omega} g(\textbf{v})d\textbf{v} = 0$ for any
function $g$ odd in $\textbf{v}$. Then $f$ solves the non dimensional linear
transport equation
\begin{eqnarray}\label{TEprinc}
\varepsilon \partial_t f + \textbf{v} \nabla_x f =
\frac{1}{\varepsilon}\left( \frac{\sigma_s}{S}\int_{\Omega} f
d\textbf{v}' - \sigma f\right) + \varepsilon Q,
\end{eqnarray}
where $\sigma = \sigma(x)$ is the total cross section, $\sigma_s =
\sigma_s(x)$ is the scattering coefficient. Here $Q = Q(x) $ is a
source term and $\varepsilon$ the mean free path. Typically,
$\sigma_s = \sigma - \varepsilon^2 \sigma_A$ where $\sigma_A =
\sigma_A(x)$ is the absorption coefficient. Such an equation
arises in neutron transport \cite{CZ}, radiative transfer 
\cite{Cha} and wave propagation in random media \cite{RPK}, etc.
In all these applications, the scaling appearing in  (\ref{TEprinc}) is
typical, and gives rise to a diffusion equation as $\varepsilon
\rightarrow 0$ of the form \cite{JPT2}
\[
\partial_t \rho=\frac1{S}\int_{\Omega} \textbf{v}\cdot \nabla_x \left(\frac{\textbf v}{\sigma}\cdot\nabla_x\rho\right)\,d\textbf{v}-\sigma_A\rho+{Q},
\]
where
$\rho=(1/S)\int_{\Omega}f\,d\textbf{v}$. We refer, for example, to \cite{BLP, BSS} and the references therein for rigorous mathematical results concerning diffusion limits of transport equations.

\subsection{Problem reformulation}
Consider now the one-dimensional transport equation
\begin{eqnarray}\label{TE}
\varepsilon \partial_t f + v \partial_x f =
\frac{1}{\varepsilon}\left( \frac{\sigma_s}{2}\int_{-1}^{1} f
d{v}' - \sigma f\right) + \varepsilon Q,
\end{eqnarray}
with $x_L< x < x_R$ and boundary conditions
\begin{eqnarray}
\begin{array}{cc}
f(t,x_L, v) = F_L(v), & \textrm{for} \ \ v > 0,\\
f(t,x_R, -v) = F_R(v), &  \textrm{for} \ \ v > 0.
\end{array}
\end{eqnarray}
In \cite{JPT2} the authors proposed a method based on the even-odd
decomposition $f = r + \varepsilon j$ where $r = \frac{1}{2} (f(v)
+ f(-v))$ and $j = \frac{1}{2 \varepsilon} (f(v) - f(-v))$, that
splits the equation (\ref{TE}) as two equations, each for $v>0$
 \begin{eqnarray}\label{ET2}
\varepsilon \partial_t f(v) + v\partial_x f(v) = \frac{1}{\varepsilon}\left( \frac{\sigma_s}{2}\int^{1}_{-1} f dv' - \sigma f(v)\right) + \varepsilon Q, \\
\varepsilon \partial_t f(-v) - v\partial_x f(-v) =
\frac{1}{\varepsilon}\left( \frac{\sigma_s}{2}\int^{1}_{-1} f dv'
- \sigma f(-v)\right) + \varepsilon Q. \nonumber
\end{eqnarray}
Adding and subtracting these two equations leads to
\begin{eqnarray}\label{equaMain}
\partial_t r + v\partial_x j = -\frac{\sigma_s}{\varepsilon^2}(r-\rho)-\sigma_A r + Q,\\ \nonumber
\acapo
\partial_t j + \frac{v}{\varepsilon^2}\partial_x r = -\frac{\sigma_s}{\varepsilon^2}j-\sigma_A j,\nonumber
\end{eqnarray}
where
\begin{eqnarray}
\rho=\int_{0}^{1} r dv.
\end{eqnarray}
 As $\varepsilon \rightarrow 0$, system (\ref{equaMain}) gives
\begin{eqnarray}\label{lim1}
r = \rho, \ \ j = -(v/\sigma) \partial_x r. \nonumber
\end{eqnarray}
Applying this to the first equation of (\ref{equaMain}) and
integrating over $v$ we get the diffusion equation
\begin{eqnarray}\label{difflim}
\partial_t\rho = \frac{1}{3} \partial_{xx} \rho-\sigma_A \rho + Q.
\end{eqnarray}
To get boundary conditions for $r$ and $j$ we use
relations
\begin{eqnarray}\label{bc1}
r + \varepsilon j|_{x = x_L} = F_L(v), \ \ \ r - \varepsilon j|_{x
= x_R} = F_R(v).
\end{eqnarray}

For arbitrary value of $\varepsilon$, the compatibility conditions imposed at
the boundary is quite complicated. However, for small values of $\varepsilon$,
which is the case we are interested in, the treatment simplifies, because 
when $\varepsilon \rightarrow 0$, $j = -(v/\sigma) \partial_x r$,
then applying this in (\ref{bc1}) one gets
\begin{eqnarray}\label{bc}
r - \varepsilon v\partial_x r|_{x = x_L} = F_L(v), \ \ \ r +
\varepsilon v \partial_x r|_{x = x_R} = F_R(v).
\end{eqnarray}

Such boundary conditions will avoid a boundary layer in the limit
case $\varepsilon \rightarrow 0$, therefore the numerical boundary
conditions are obtained by discretizing Eq. (\ref{bc}). As done in
\cite{JPT2}, the boundary conditions have been applied using a
second order implementation of equation (\ref{bc1}) based on
central differences. Extensions to higher-order
implementation of equations (\ref{bc}) and to different boundary conditions 
are not considered here and will be investigated in a forthcoming work.
According to this, for our tests we chose a second order scheme in space and time, 
because the boundary conditions are discretized
to second order accuracy.

Now we start from system (\ref{equaMain})
and, adding and subtracting the quantity $v^2\partial_{xx} r/\sigma$
in the first equation, we reformulate the problem in the equivalent form
\begin{eqnarray}\label{equaMain3}
\begin{array}{ll}
\partial_t r = \underbrace{- v\partial_x( j + \frac{\mu(\varepsilon)v\partial_x r}{\sigma})}_{Explicit} \underbrace{-\frac{\sigma_s}{\varepsilon^2}(r-\rho)-\sigma_A r + Q+\mu(\varepsilon)v^2\frac{\partial_{xx}r}{\sigma}}_{Implicit},\\
\partial_t j  = \underbrace{- \frac{1}{\varepsilon^2}(j + \frac{v\partial_x r}{\sigma})}_{Implicit},% -\frac{\sigma_s}{\varepsilon^2}j-\sigma_A j,
\end{array}
\end{eqnarray}
The new system is then discretized in time using an IMEX R-K
scheme as described in the first part of the manuscript. The implicit-explicit integration process is emphasized in (\ref{equaMain3}).

We remark that this new formulation of the diffusive relaxation
system (\ref{equaMain}) is such that when $\varepsilon$ tends to
zero the system relaxes towards (\ref{lim1}). From a numerical
point of view the new formulation has several advantages. In particular, 
as $\varepsilon \rightarrow 0$, the IMEX R-K scheme applied to 
system (\ref{equaMain3}) originates a
fully implicit scheme for solving the diffusion equation
(\ref{difflim}).

%From now on, we will denote this new reformulation \emph{BPR approach} in the numerical tests. In this approach, 
In our numerical tests  in order to obtain uniformly accurate second order scheme both in space and 
in time for the BPR approach, we consider SSP2(3,3,2)-WENO32. We remark that we can 
obtain analogous results considering ARS(2,2,2)-WENO32.

The equations are discretized in space and velocity, i.e.
$r(x_i,v_m,t_n) \approx r^{n}_{i,m}$ where $\left\{v_m\right\}$
are chosen to be the $N_v$ positive nodes of the Gauss-Legendre
quadrature formula, with $2N_v$ nodes in the interval $[-1, 1]$
while $x_i = \Delta x(i-1/2)$ for $i = 1,...,N_p$. Note that the
computation of the $k$-th stage of the implicit equation
$r^{(k)}_{i,m}$, requires the quantity $\rho^{(k)}_{i}$ in the
implicit part in (\ref{equaMain3}). Such quantities are obtained
as follows. Assume we have computed $r^{(l)}_{i,m}$ for $l =
1,...,k-1$, then $r^{(k)}_{i,m}$ is obtained from
\begin{eqnarray}\label{equref}
r^{(k)}_{i,m} = \overline{r}^{(k-1)}_{i,m} + \Delta t
a_{kk}\left(\frac{\sigma_S}{\varepsilon^2}(r^{(k)}_{i,m} -
\rho^{(k)}_i)-\sigma_A r^{(k)}_{i,m} + Q\right)
\end{eqnarray}
discretizing system  (\ref{equaMain3}) without the quantity
$\mu(\varepsilon)v^2\partial_{xx} r/\sigma$ in the implicit and explicit part. The
quantity $\overline{r}^{(k-1)}_{i,m}$ represents the contribution of the first $k-1$ stages. 
Then, in order to compute $\rho^{(k)}_i$ we apply Gauss
quadrature on both sides of (\ref{equref}) (i.e. multiply by the
weights $w_m$ and sum over $m$), setting $\mu = 0$. In this way we
obtain an equation for $\rho^{(k)}_i$ that can be explicitly
solved, and such value is plugged in (\ref{equaMain3}) in order to
compute $r^{(k)}_{i,m}$.

%%%%%%%%%%%%%%%%%%%%%%%%%%%%%%%ESEMPI%%%%%%%%%%%%%%%%%%%%%%%%%%%%%
\subsection{Numerical results}
In this section we shall consider some transport problems in slab
geometry. We will present the transient and the steady state
solutions. We remark that in all the test problems we have used
$N_v=8$ thus the standard $16$ points Gaussian quadrature set for the velocity space.
In all the tests the initial distribution is $f(\textbf{x},\textbf{v},t=0)=0$.

We emphasize that, besides uniform accuracy in $\varepsilon$, our approach 
allows to choose larger time steps, since there is no stability restriction on the time step.
As we will show, this permits to obtain numerical results at a lower computational cost
compared to other approaches presented in the literature that lead to explicit schemes for the underlying diffusion limit
with a parabolic CFL stability restriction \cite{JPT2, Klar}.  
Nevertheless, in order to get an accurate resolution of the behavior of the solution, smaller time step may
be necessary.

Depending on the regime of the parameter $\varepsilon$, we compare numerical solutions
to a direct implicit discretization of the diffusion limit
(\ref{difflim}) when $\varepsilon$ tends to zero, whereas for intermediate values of the parameter
$\varepsilon$  we compute a reference solution using a much finer grid in space. 

In the next tests we compare the results obtained by the new approach versus the results given by   
Jin et al., in \cite{JPT2}, here denoted by JPT. We refer to \cite{JPT2, Klar, LM} for similar results where the limiting scheme is explicit and so in diffusive regions requires $\Delta t \approx (\Delta x)^2$. In all
figures we use notations $N_s$ and $N_p$ to denote the number of time steps and grid points in space respectively.

\begin{figure}%[ht]
\centering
%\vspace{2.5in}
\includegraphics[width=0.49\textwidth]{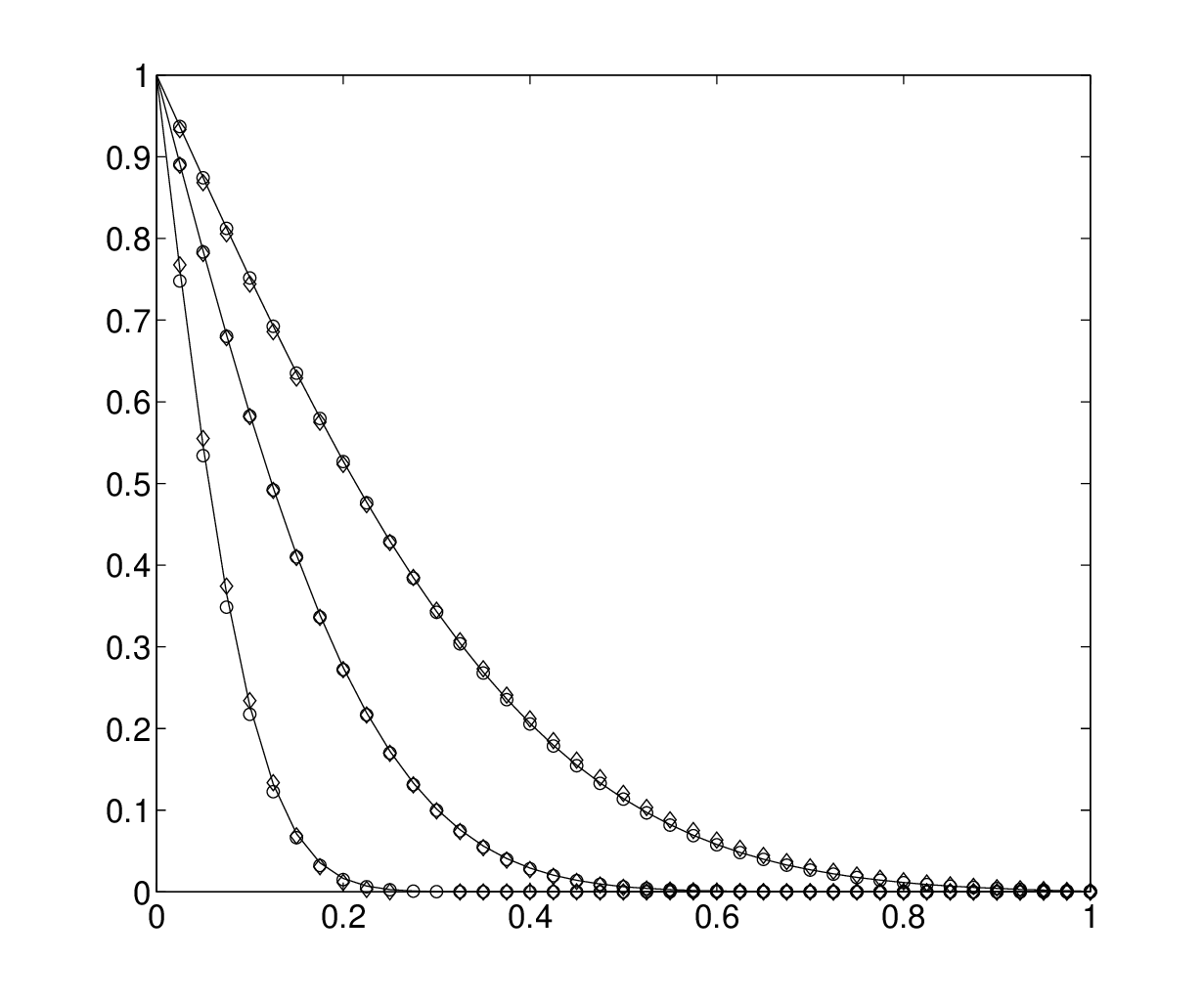}
\includegraphics[width=0.5\textwidth]{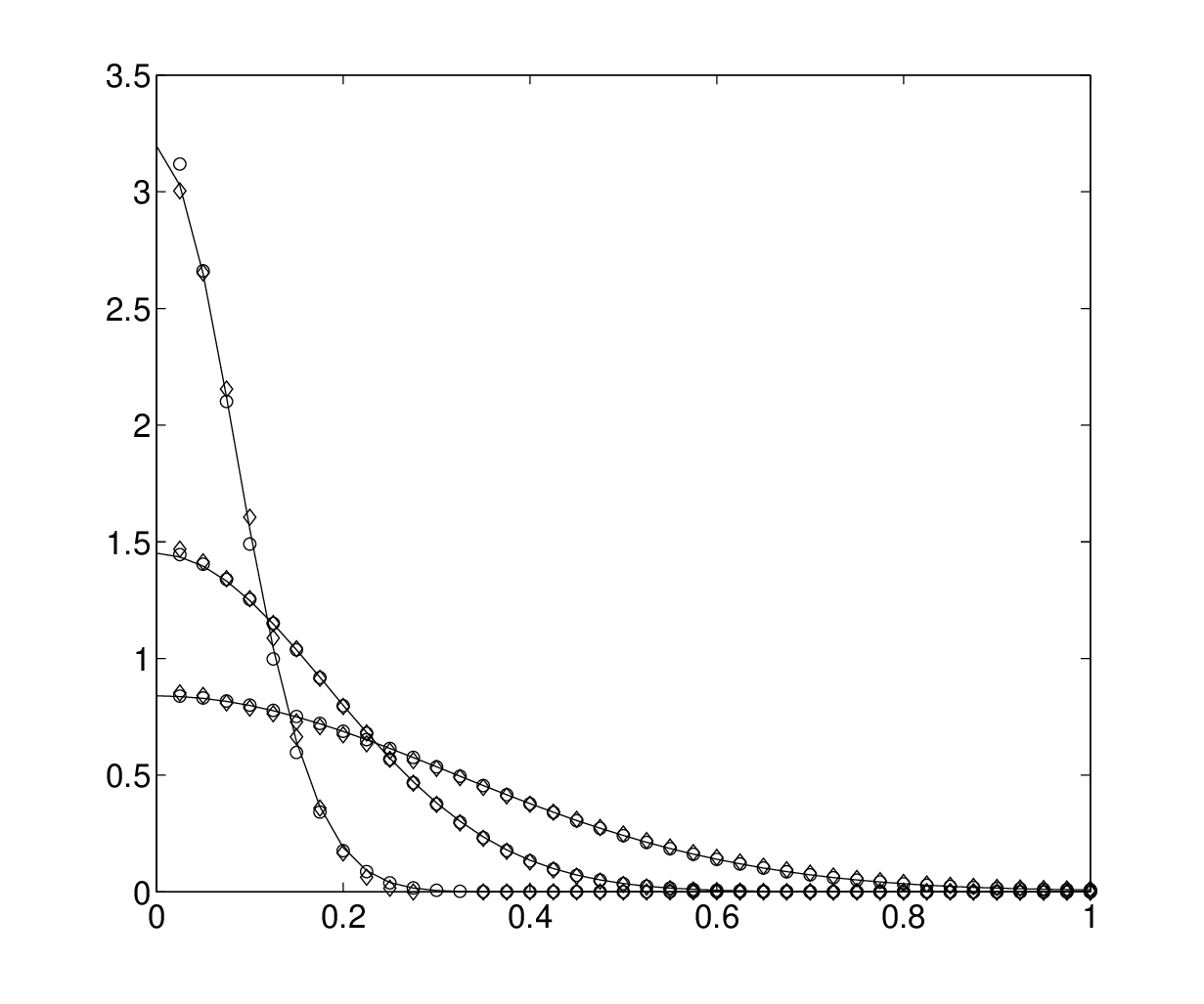}
\caption{Problem I. On the left-hand side  the mass density $\rho$, on the right-hand side the flux $\rho u$.
JPT ($\circ$), $\Delta x = 0.025$,
$\Delta t = 0.0002$, $N_s = 50, \ \  250, \ \ 750$.  BPR
($\Diamond$),  $\Delta t = \lambda \Delta x$, with $\lambda = 0.035$,
$N_s = 11, \ \  57, \  \ 171$. The exact diffusive solution is
represented by the solid line.}
\label{fig:4}
\end{figure}

\begin{figure}[!ht]
\centering
\caption{The steady state solution for the Problem I and II.} \subfigure[Problem I. The steady state solution of the mass density $\rho$.  JPT ($\circ$),  $\Delta t = 0.001$, $N_s = 2000$.  BPR ($\Diamond$), $\Delta t = 0.0035$, $N_s = 561$. The exact diffusive solution is represented by the solid line.]{
   \includegraphics[width=0.45\textwidth]{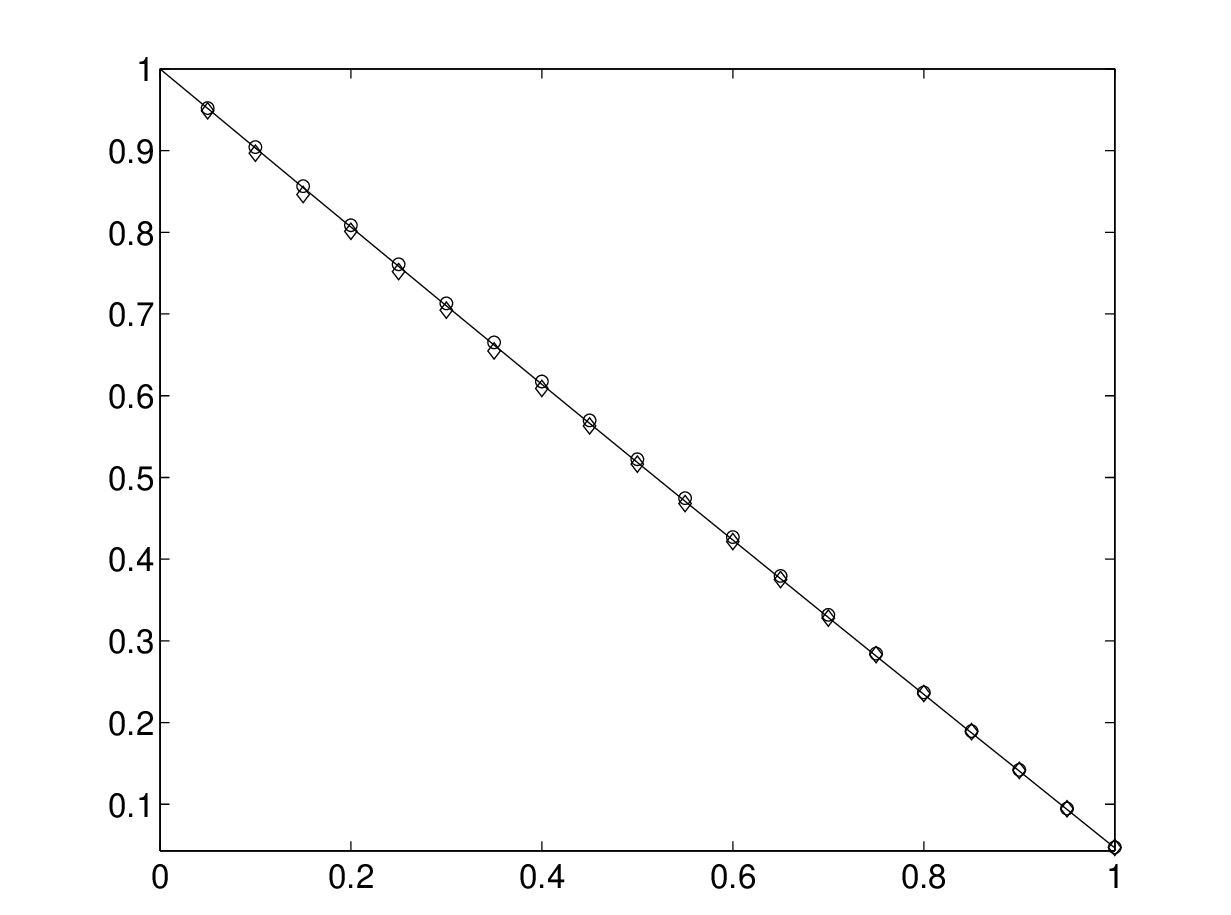}}\label{fig:5}
 \hspace{7mm}
 \subfigure[Problem II. The steady state solution of the mass density $\rho$ with $\varepsilon = 1$ and $\Delta x = 0.05$  on $(0, 1)$ and  with $\varepsilon = 0.01$ and $\Delta x = 1$ on $(1, 11)$. JPT ($\circ$), with $\Delta t = 0.025 $ and $N_s = 6000$. BPR ($\Diamond$), with $\Delta t = \Delta x$, with $N_s = 3000$]{
   \includegraphics[width=0.45\textwidth]{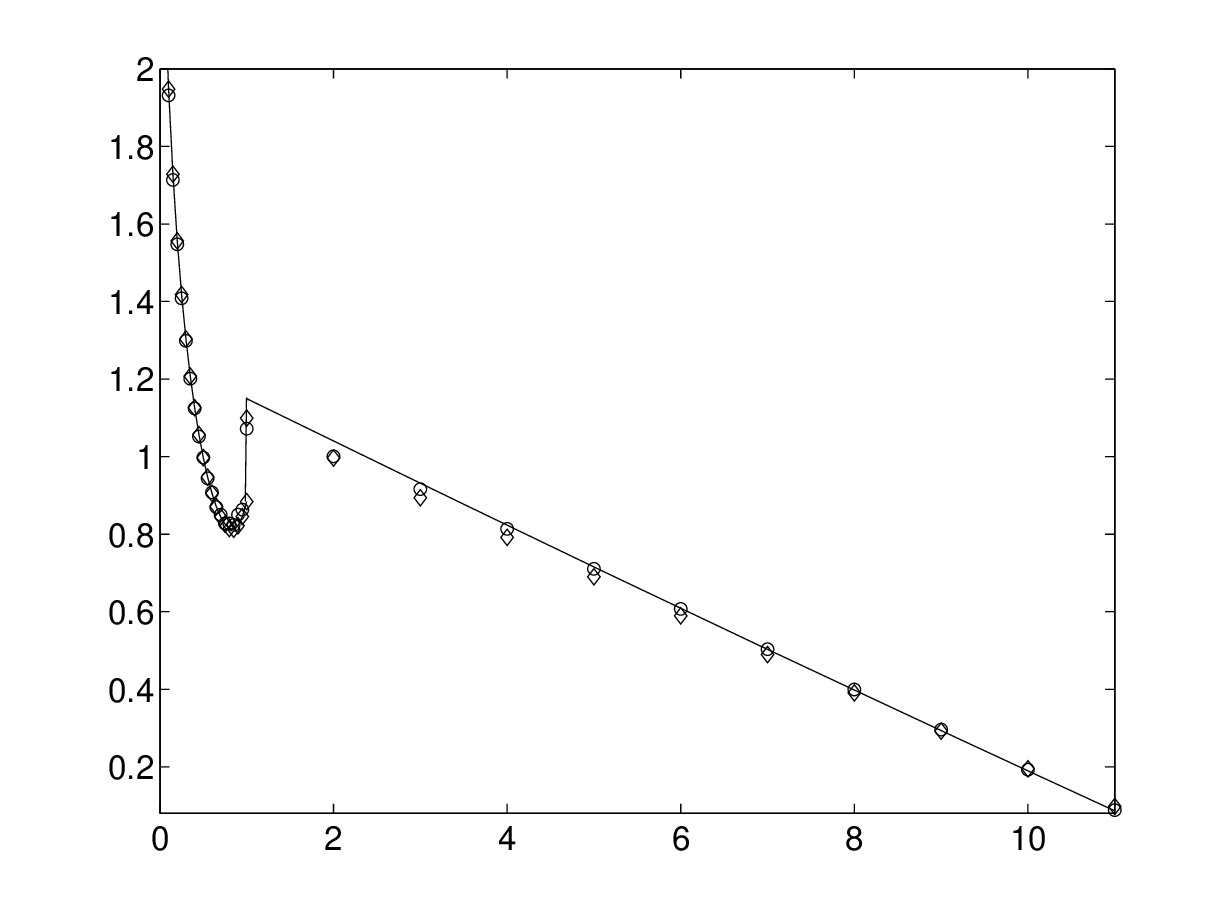}}\label{fig:6}
\end{figure}

\paragraph{Problem I:}
\begin{eqnarray*}
\begin{array}{ccc}
x \in [0, 1], & F_L(v) = 1, & F_R(v) = 0,\\
\sigma_S = 1, \ \ \ \sigma_A = 0, & Q = 0, & \varepsilon =
10^{-8}.
\end{array}
\end{eqnarray*}
The numerical results are reported in Figures \ref{fig:4} and \ref{fig:5} (a) at different
times $t = 0.01, \ 0.05, \ 0.15$ with $N_p = 40$, and at $ t = 2$ with $N_p = 20$ where the steady
state is reached. In this problem we see that in both cases, the
results in the transient and steady state solutions show a good
behavior with the correct diffusion limit. The exact diffusive solution has been computed
by  (\ref{difflim}) with $N_p = 200$. As expected both
JPT and BPR results are very close to the exact diffusive solution at any
times. Note however that thanks to the better stability properties in this regime BPR scheme is about $4$ times faster then the explicit method.
%%%%%%%%%%%%%%%
\paragraph{Problem II:} This is a two-material problem used in
\cite{JPT2, LM, Klar} where in the purely absorbing region  $[0,
1]$ the solution decays exponentially whereas in the purely
scattering region  $[1, 11]$ the solution is diffusive, the
parameters are the following
\begin{eqnarray*}
x \in [0, 11], \ \ \ F_L(v) = 5, \ \ \ F_R(v) = 0,\\
\sigma_S = 0, \ \ \ \sigma_A = 1, \ \ \ Q = 0, \ \ \varepsilon = 1, \ \  \textrm{for} \ \ x \in [0,1],\\
\sigma_S = 1, \ \ \ \sigma_A = 0, \ \ \ Q = 0, \ \ \varepsilon = 0.01, \ \ \textrm{for} \ \ x \in [1,11].\\
\end{eqnarray*}
An interface layer is produced between the pure absorbing region
and the scattering one. Two meshes are used in the domain $[0,
11]$, a thin mesh $\Delta x = 0.05$ in $[0, 1]$ and a coarse mesh
$\Delta x = 1$ in $[1, 11]$, which means that between the
interface layer we have to use a space discretization with a non
uniform mesh. Finite volume, rather than finite difference, is used in this case.
Since we restrict to second order accuracy, the point-wise values of the source is identified  with its cell average.
 The high order non oscillatory reconstruction is performed
by a WENO approach for non uniform mesh. For our numerical tests 
we used WENO32 reconstruction. 

\begin{figure}%[ht]
\centering
%\vspace{2.5in}
\includegraphics[width=0.5\textwidth]{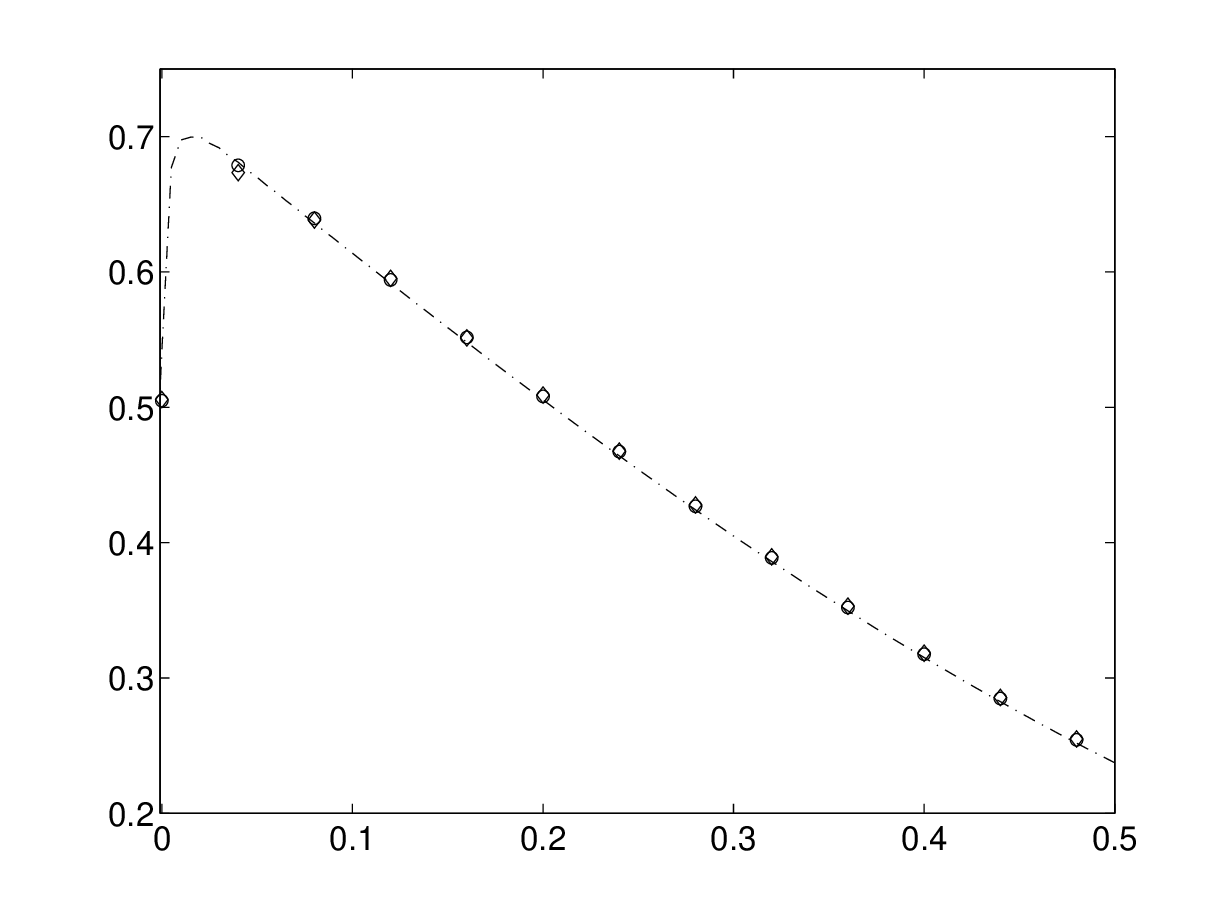}%{ProblemaIVrho2.eps}
\includegraphics[width=0.5\textwidth]{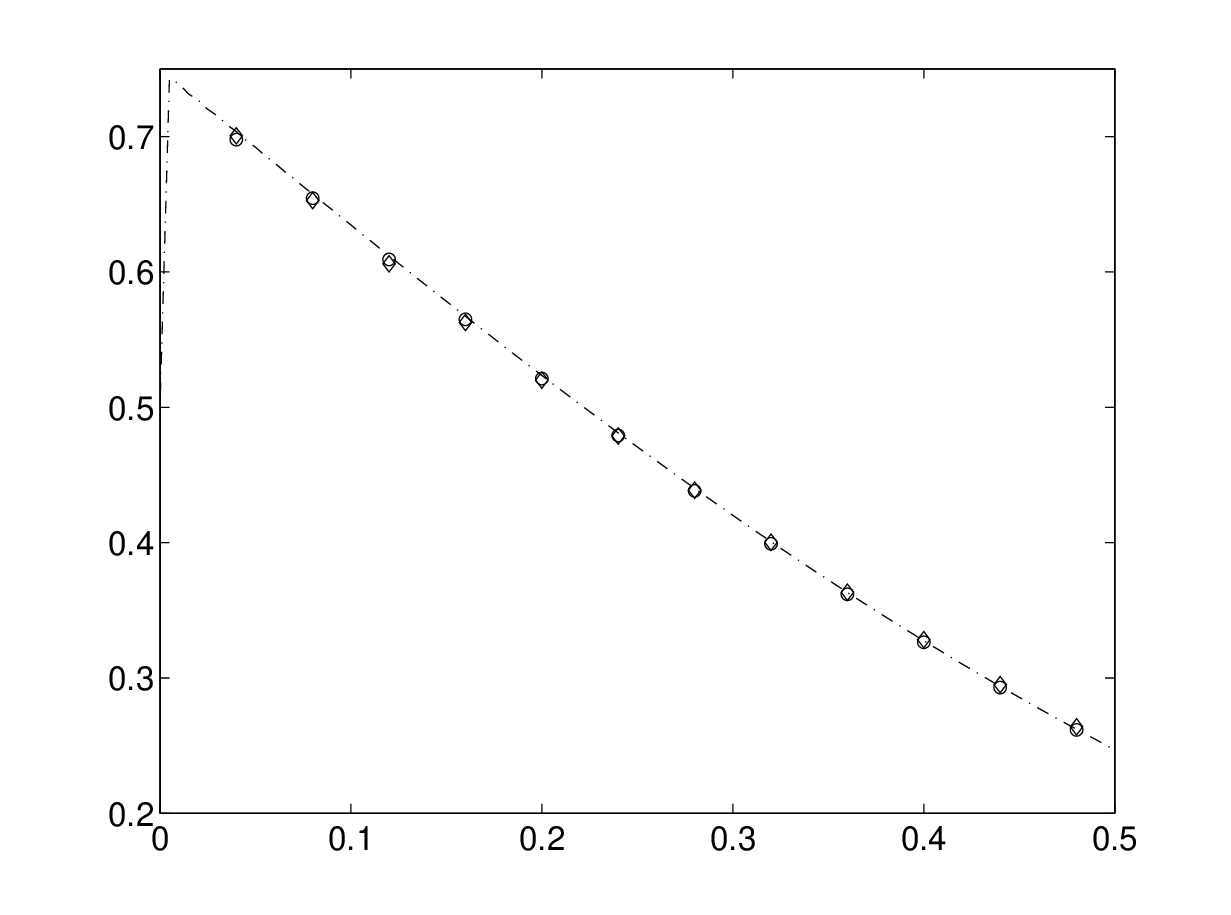}%{ProblemaIVrho1e_42.eps}
\caption{Problem III. The mass density $\rho$. JPT ($\circ$), with $\varepsilon = 10^{-2}$, $\Delta x = 0.04$, $\Delta t = 0.001$, $N_p = 25$, $N_s = 400$. BPR ($\Diamond$), $\Delta x = 0.04$, $\Delta t= 0.002 $, $N_p = 25$, $N_s = 100$. Left $\varepsilon = 10^{-2}$, right $\varepsilon=10^{-4}$. The reference solutions are represented by the dash dot line.}
\label{fig:7}
\end{figure}

At time $t = 150$ the solution has reached the steady state and the results are presented in Figure \ref{fig:6}. We
computed the reference solution with a very fine discretization $N_p = 400$
using a uniform mesh in all the domain $[0, 11]$. The numerical
schemes provide a good description for the solution in the
absorption and diffusive regions, in fact, we observe that for the
steady state, Figure \ref{fig:6} (b), JPT and BPR results are
close o the reference solution, except at the interface where
a slight difference is observed. In this case BPR scheme is about twice times faster then JPT approach.
%%%%%%%%%%%%%%%%%%%%
%%%%%%%%%%%%%%%%%%%%%%%%
\paragraph{Problem III:} Concerning this problem we present two
different situations (see \cite{JPT2, LM}) with non-isotropic boundary conditions that generate a
boundary layer:
\begin{eqnarray*}
x \in [0, 1], \ \ \ F_L(v) = v, \ \ \ F_R(v) = 0,\\
\sigma_S = 1, \ \ \ \sigma_A = 0, \ \ \ Q = 0, \ \ \ \varepsilon =
10^{-2}.
\end{eqnarray*}
First in an intermediate regime with $\varepsilon=10^{-2}$ and then in a more diffusive regime with $\varepsilon= 10^{-4}$.
Using a coarse discretization $N_p = 25$ the boundary layer is not
resolved, but we observe that the two approaches accurately capture the solution inside the domain
(in Figure (\ref{fig:7}) we have restricted the numerical and the reference solution to the interval $[0, 0.5]$). 
The reference solution has been
obtained using a fine discretization $N_p= 400$ and the boundary layer is
resolved. The results are plotted at time $t = 0.4$ in figure \ref{fig:7}. The higher efficiency of the present method results in an improved time step ratio of a factor $4$. 

%%%%%%%%%%%%%%%%%%%%%%%%%%%%%%%
%%%%%%%%%%%%
%%%%%%%%%%%%%
%%%%%%%%%%%%

\section{Conclusions}\label{Sec7}
In this manuscript we have presented a general way to tackle diffusion limit for hyperbolic and kinetic problems which permits to obtain accurate and efficient schemes both in rarefied and diffusive regimes. The new approach, in particular, give rise to a fully implicit method for the diffusion component of the limiting system. This is obtained without solving nonlinear systems of implicit equations but by a suitable blending into the IMEX R-K method of a fully implicit solver for the limiting diffusive system. Numerical results show that this approach is able to capture the correct asymptotic behavior of the system at a lower computational cost compared to other approaches that lead to explicit schemes for the underlying diffusion limit, because we removed the parabolic CFL restriction, common to most approaches in the literature.

The method here presented is based on the use of IMEX Runge-Kutta methods, however extension to more general additive Runge-Kutta schemes are naturally possible. In fact from our problem reformulation
\begin{eqnarray}\label{SPPc}
\begin{array}{l}
\displaystyle u' = \underbrace{f_1(u,v)}_{explicit} + \underbrace{f_2(u)}_{implicit}, \\ [+.25cm]
\displaystyle {\varepsilon^2}v' =  \underbrace{g(u,v)}_{implicit}\\
\end{array}
\end{eqnarray}
we can clearly combine different implicit solvers to tackle the ``highly" stiff component $g(u,v)$ which originates the algebraic condition $g(u,v)=0$ giving rise to the equilibrium projection $v=G(u)$ and the ``mildly" stiff component $f_2(u)$ corresponding to the limiting diffusive term. We leave this possibility to future research directions.

\section{\large{Appendix}}

\subsection{Stability analysis of first order IMEX schemes}\label{A1}

For the subsequent analysis we restrict to the linear case
$p(u)=u$ and $q(u)=0$ \begin{eqnarray}
\begin{array}{l}
\displaystyle
u_t  =- (v+\mu u_x)_x + \mu u_{xx}, \\[+.25cm]
\epsi^2\displaystyle v_t =- u_x -v,
\\
\end{array}
\label{I3c} \end{eqnarray}

We now look for a Fourier solution of the form $u = \hat
u(t)\exp(i\xi x)$, $v = \hat v(t)\exp(i\xi x)$ and inserting the
{\em ansatz} into systems \rf{I3}, the evolution equations are
\begin{eqnarray}
  \hat u_t & = & - i \xi \hat v + \xi^2 \mu \hat u -\xi^2 \mu \hat u    \nonumber         \acapo
            &   &                    \label{eq:relaxF2} \acapo
 \varepsilon^2 \hat v_t & = & - i \xi \hat u  -\hat v      \nonumber
\end{eqnarray}
It is convenient to rewrite the system using the variable $\hat w  = {-i\hat v}/{\xi}$ in place of $\hat{v}$ so the system becomes
\begin{eqnarray} \label{finalversion}
  \hat u_t & = & \xi^2( \hat w +  \mu \hat u) -\xi^2 \mu \hat u,    \nonumber     \acapo
            &   &                    \label{eq:relaxF3} \acapo
 \varepsilon^2 \hat w_t & = & - \hat u  -\hat w.      \nonumber
\end{eqnarray}

We apply the first order IMEX method based on explicit and
implicit Euler schemes to get
\begin{eqnarray} \nonumber
    \hat u^{n+1} &=&\hat u^n + {\Delta t}\xi^2( \hat w^n +  \mu\hat u^n) -{\Delta t}\xi^2 \mu \hat u^{n+1},    \nonumber     \acapo
            &   &                    \label{IEE} \acapo
 \varepsilon^2 \hat w^{n+1}&=&\varepsilon^2  \hat w^n - {\Delta t}\hat u^{n+1}  -{\Delta t}\hat w^{n+1},      \nonumber
\end{eqnarray}
which after manipulation can be written explicitly in the form
\begin{eqnarray} \nonumber
    \hat u^{n+1} &=&\hat u^n + \frac{\Delta t\xi^2}{1+\Delta t\xi^2 \mu}\hat w^n,    \nonumber     \acapo
            &   &                    \label{IEE2} \acapo
            \nonumber \acapo
\hat w^{n+1}&=&\frac{\varepsilon^2(1+\Delta t\xi^2 \mu)-\Delta
t^2\xi^2}{(\varepsilon^2+\Delta t)(1+\Delta t\xi^2 \mu)}\hat w^n -
\frac{\Delta t}{\varepsilon^2+\Delta t}\hat u^{n}. \nonumber
\end{eqnarray}
In order to study the stability of the method we compute the
eigenvalues of the stability matrix
\begin{equation}
R=\left(%
\begin{array}{cc}
  1 & \displaystyle\frac{\Delta t\xi^2}{1+\Delta t\xi^2 \mu}
  \\[+.5cm]
  -\displaystyle\frac{\Delta t}{\varepsilon^2+\Delta t} & \displaystyle\frac{\varepsilon^2(1+\Delta t\xi^2 \mu)-\Delta
t^2\xi^2}{(\varepsilon^2+\Delta t)(1+\Delta t\xi^2 \mu)} \\
\end{array}%
\right).
\end{equation}
We obtain the expressions
\begin{equation}
\lambda_{\pm}=\frac12\left\{1+\alpha(1+\beta)-\beta\pm\sqrt{[1+\alpha(1+\beta)-\beta]^2-4\alpha}\right\}
\end{equation}
with
\[
\alpha=\frac{\varepsilon^2}{\varepsilon^2+\Delta t},\quad
\beta=\frac{\Delta t\xi^2}{1+\mu\Delta t\xi^2}.
\]
It can be shown that $|\lambda_{\pm}|<1$ when
\[
\alpha < \frac{1-2\beta+\beta^2}{1+2\beta+\beta^2}.
\]
The above inequality involves a third order polynomial in $\Delta
t$
\begin{eqnarray}
 -\xi^4(\mu-1)^2{{\Delta t
}}^{3}+ 2\xi^2(2\,{\varepsilon}^{2}{\xi}^{2}\mu+1-\mu) {{\Delta
t}}^{2}+ \left( 4\,{\varepsilon}^{2}{\xi}^{2}-1
 \right) {\Delta t}<0.
 \label{cond1}
\end{eqnarray}
The roots of this polynomial are given by $T_0=0$ and
\[
T_{\pm}=\frac1{(\mu-1)^2\xi^2}\left\{2\varepsilon^2\xi^2\mu-\mu+1\pm
2\varepsilon\xi\sqrt{\varepsilon^2\xi^2\mu^2-\mu+1}\right\}.
\]
Condition (\ref{cond1}) can be satisfied only if the last two
roots are positive. This is guaranteed when $2\varepsilon|\xi|<1$ and
so we have the time step restriction $\Delta t < T_{-}$. The
largest stability region is obtained when $\mu=1$ for which we get
\begin{equation}
\xi^2\Delta t < \frac14
\frac{(1-4\xi^2\varepsilon^2)}{\varepsilon^2\xi^2},
\end{equation}

This relation is interpreted as follows. For a fixed $\varepsilon |\xi|< 1/2$, the restriction on the time 
step is of parabolic type, since $|\xi| \sim 1/\Delta x$ is the maximum Fourier mode represented on a grid of spacing $\Delta x$.
The restriction on $\Delta t/\Delta x^2$ is less and less severe as $\varepsilon \rightarrow 0$.

The implicit treatment of the second equation stabilizes the explicit treatment of the first one, 
provided $\varepsilon$ is sufficiently small.  
%where the right hand side is a decreasing function of
%$\varepsilon^2$.

%%%%%%%%%%%%%%%%%%%%%%%%%%%%5
\subsection{Analysis of second order stiffly accurate schemes}\label{A2}
%%%%%%%%%%%%%%%%%%%%%%%%%%%%5
%%%%%%%%%%%%%%%%%%%%%%%%%%%%%
We have the following result:
\begin{theorem}\label{prop}
Consider an IMEX Runge-Kutta scheme of type A. Then
there exist no second-order tree stage scheme satisfying the
conditions $\bt A^{-1} = e_s^T$ and $\tbt = e_s^T\tA$.
\end{theorem}\\
{\bf Proof.} We consider the classical second order conditions
\begin{eqnarray}\label{secondorder}
\begin{array}{cc}
\tbt e = 1, & \bt e = 1,\\
\tbt \tc = 1/2, & \bt c = 1/2,\\
\tbt c = 1/2, & \bt \tc = 1/2,
\end{array}
\end{eqnarray}
with $c = A \one$ and $\tc = \tA \one$ and the conditions $\bt A^{-1} = e_s^T$ and $\tbt = e_s^T\tA$.\\
For $s = 3$ the Butcher \emph{tableau} of a stiffly accurate IMEX R-K of type A is
\begin{displaymath}
\begin{array}{c|ccc}
 0&0 & 0& 0\\
\tilde{c}_2 &  \tilde{c}_2  & 0 & 0 \\
1 & \tilde{b}_1&\tilde{b}_2 & 0\\
\hline
 &  \tilde{b}_1 & \tilde{b}_2 &0\\
\end{array}
\quad
\begin{array}{c|ccc}
 c_1 & c_1 & 0 & 0\\
 c_2&  c_2 - a_{22} & a_{22} & 0          \\
  1 &  b_1 & b_2 & \gamma \\
\hline
 & b_1 & b_2 & \gamma
\end{array}
\end{displaymath}
(note that stiff accuracy implies $c_s = \tilde{c}_s = 1$)
and from (\ref{secondorder}) the resulting system of equations can be explicitly written
\begin{eqnarray}\label{solve2}
\begin{array}{cc}
\tb_1  = 1-\tb_2, & b_1  = 1-\gamma - b_2,\\
\tb_2 \tc_2 = 1/2, & (1- \gamma- b_2)c_1 + b_2 c_2 = 1/2-\gamma,\\
(1-\tb_2)c_1 + \tb_2 c_2 = 1/2, &  b_2 \tc_2 = 1/2-\gamma.
\end{array}
\end{eqnarray}
In order to solve system (\ref{solve2}) we can compute the coefficients as follows
\begin{eqnarray*}
\tb_2  =1/(2\tc_2), \quad b_2 = (1-2\gamma)/(2\tc_2),
\end{eqnarray*}
and
\begin{eqnarray}\label{solve3}
\begin{array}{cc}
b_2(c_2-c_1) = 1/2 -\gamma -c_1 +c_1\gamma,\\
\tb_2(c_2 - c_1) = 1/2 -c_1.
\end{array}
\end{eqnarray}
Substituting $\tb_2$ and $b_2$ in (\ref{solve3}) we get
\begin{eqnarray*}
\frac{(c_2-c_1)}{2 \tc_2} &=& \frac{1/2 -\gamma -c_1 +c_1\gamma}{1-2 \gamma},\\
\frac{(c_2 - c_1)}{2 \tc_2} &=& 1/2 -c_1.
\end{eqnarray*}
Now, comparing and equating the two expressions we have $\gamma= 0$ and it is impossible because the matrix $A$ is invertible.\\
\subsection{Second and third order IMEX schemes}\label{A4}
\begin{enumerate}
\item Second order IMEX schemes:
\begin{itemize}
\item[-] ARS(2,2,2) scheme, \cite{ARS}
\begin{eqnarray*}
\begin{array}{c|ccc}\footnotesize
0 & 0 &0 & 0\\
\gamma & \gamma &0& 0\\
1 & \delta & 1-\delta & 0\\
\hline
 & \delta & 1-\delta & 0
\end{array} \ \ \ \ \
\begin{array}{c|ccc}
0 & 0 & 0 & 0\\
\gamma & 0 &\gamma & 0\\
1 & 0 & 1-\gamma & \gamma\\
\hline
 & 0 & 1-\gamma & \gamma
\end{array}
\end{eqnarray*}
with 	$\gamma = (2-\sqrt{2})/2$ and $\delta =1 -1/(2 \gamma)$.

\item[-] SSP2-(3,3,2) scheme, \cite{A}
\begin{eqnarray*}
\begin{array}{c|ccc}\footnotesize
0 & 0 & 0& 0\\
1/2 & 1/2 & 0& 0\\
1 & 1/2 & 1/2 & 0\\
\hline
 & 1/3 &1/3 & 1/3
\end{array} \ \ \ \ 
\begin{array}{c|ccc}
1/4 & 1/4 & 0 & 0\\
1/4 & 0 & 1/4& 0\\
1 & 1/3 & 1/3 & 1/3\\
\hline
 & 1/3 &1/3 & 1/3
\end{array}
\end{eqnarray*}
\end{itemize}

\item Third order IMEX schemes:
\begin{itemize}
\item[-] ARS(4,4,3) scheme, \cite{ARS}
\begin{eqnarray*}
\begin{array}{c|ccccc}\footnotesize
0 & 0 &0 &0 & 0& 0\\
1/2 & 1/2 &0 &0 & 0 & 0\\
2/3 & 11/18 & 1/18 &  0&0 &0\\
1/2 & 5/6 & -5/6 & 1/2 & 0&0\\
1 & 1/4 & 7/4 &  3/4 & -7/4&0\\
\hline
 & 1/4 & 7/4 & 3/4 & -7/4 &0
\end{array} \ \ \ \ \
\begin{array}{c|ccccc}
0 & 0 & 0 & 0 & 0& 0\\
1/2 & 0& 1/2 & 0 & 0\\
2/3  &0 & 1/6 & 1/2 &0 &0\\
1/2 & 0 & -1/2 & 1/2&1/2 & 0\\
1  & 0 & 3/2 & -3/2 & 1/2 &1/2\\

\hline
 & 0 & 3/2 & -3/2 & 1/2 &1/2
\end{array}
\end{eqnarray*}

\item[-] BPR(3,5,3) scheme
\begin{eqnarray*}
\begin{array}{c|ccccc}\footnotesize
0 & 0 &0 &0 & 0& 0\\
1 & 1 &0 &0 & 0 & 0\\
2/3 & 4/9 & 2/9 &  0&0 &0\\
1 & 1/4 & 0 &  3/4 & 0&0\\
1 & 1/4 & 0 &  3/4 & 0&0\\
\hline
 & 1/4 & 0 & 3/4 & 0 &0
\end{array} \ \ \ \ \
\begin{array}{c|ccccc}
0 & 0 & 0 & 0 & 0& 0\\
1 & 1/2 & 1/2 & 0 & 0\\
2/3 & 5/18 & -1/9 & 1/2&0 & 0\\
1   & 1/2 & 0 & 0 &1/2 & 0\\
 1 &1/4 & 0 & 3/4 &-1/2 &1/2\\
\hline
 & 1/4 & 0 & 3/4 & -1/2 & 1/2
\end{array}
\end{eqnarray*}
\end{itemize}
\end{enumerate}

\end{document}